\documentclass[review]{elsarticle}

\usepackage{color,soul}
\usepackage{amsmath}
\usepackage{amssymb}
\usepackage{mathtools}
\usepackage{algorithmic}
\usepackage{algorithm}
\usepackage{t1enc,dfadobe}

\usepackage{amssymb}
\usepackage{url}       

\begin{document}
\title{A Unified Definition and Computation \\ of Laplacian Spectral Distances}
\author{Giuseppe Patan\'e} \ead{patane@ge.imati.cnr.it}
\address{Consiglio Nazionale delle Ricerche\\
Istituto di Matematica Applicata e Tecnologie Informatiche\\
Genova, Italy}
\begin{frontmatter}
%
%
\begin{abstract}
Laplacian spectral kernels and distances (e.g., biharmonic, heat diffusion, wave kernel distances) are easily defined through a filtering of the Laplacian eigenpairs. They  play a central role in several applications, such as dimensionality reduction with spectral embeddings, diffusion geometry, image smoothing, geometric characterisations and embeddings of graphs. Extending the results recently derived in the discrete setting~\citep{PATANE-STAR2016,PATANE-CGF2017} to the continuous case, we propose a novel definition of the Laplacian spectral kernels and distances, whose approximation requires the solution of a set of inhomogeneous Laplace equations. Their discrete counterparts are equivalent to a set of sparse, symmetric, and well-conditioned linear systems, which are efficiently solved with iterative methods. Finally, we discuss the optimality of the Laplacian spectrum for the approximation of the spectral kernels, the relation between the spectral and Green kernels, and the stability of the spectral distances with respect to the evaluation of the Laplacian spectrum and to multiple Laplacian eigenvalues.
\end{abstract}
\end{frontmatter}
\section{Introduction}
Spectral kernels and distances are easily defined through a filtering of the Laplacian eigenpairs and include random walks~\citep{SINHA2013},  biharmonic~\citep{LIPMAN2010,RUSTAMOV2011}, heat diffusion~\citep{BRONSTEIN2009,BRONSTEIN2010-TOG,GEBAL2009,LAFON2006}, wave~\citep{BRONSTEIN-PAMI2011} kernels and distances. Laplacian spectral distances have been applied to shape segmentation~\citep{GOES2008} and comparison~\citep{BRONSTEIN2010-TOG,GEBAL2009,OVSJANKOV2010,SUN2009} with multi-scale and isometry-invariant signatures~\citep{Memoli05,MEMOLI2011,MAHMOUDI2009}. Main applications of the heat kernel include dimensionality reduction~\citep{BERARD1984,XIAO2010,ELGHAWALBY2015,BAHONAR2018}, diffusion geometry~\citep{BELKIN2003,SINGER2006}, graphs' embeddings~\citep{LUO2003} and analysis~\citep{GHOSH2018}, shape comparison~\citep{BOSCAINI2016,COSMO2016,RODOLA2016}, data representation~\citep{ZHU2003} and classification~\citep{NG2001,SHI2000,SPIELMAN1996}.

\subsection*{Overview and contribution}
Extending the results recently derived in the discrete setting~\citep{PATANE-STAR2016,PATANE-CGF2017} to the continuous case, we propose a novel definition of the Laplacian spectral kernels and distances, whose approximation requires the solution of a set of inhomogeneous Laplace equations. Their discrete counterparts are equivalent to a set of sparse linear systems, which are efficiently solved with iterative methods. To this end, the Laplacian spectral distances \mbox{$d^{2}(\mathbf{p}_{i},\mathbf{p}_{j}):=\sum_{n=0}^{+\infty}\vert \phi_{n}(\mathbf{p}_{i})-\phi_{n}(\mathbf{p}_{j})\vert^{2}/\rho^{2}(\lambda_{n})$} are defined by filtering the Laplacian spectrum \mbox{$(\lambda_{n},\phi_{n})_{n=0}^{+\infty}$}, where \mbox{$\rho:\mathbb{R}^{+}\rightarrow\mathbb{R}$} is a proper filter function. 

The proposed spectral kernels and distances are independent of the data dimensionality (e.g., surfaces, volumes,~$n$-dimensional data), of the discretisation of the input domain, and of the Laplace-Beltrami operator. The Laplacian spectral distances generalise the bi-harmonic~\citep{LIPMAN2010} (\mbox{$\rho(s):=s^{2}$}), wave~\citep{BRONSTEIN-PAMI2011} (\mbox{$\rho_{t}(s):=\exp(ist)$}), and diffusion~\citep{CRANE2013,HAMMOND2010,PATANE-CAGD2013} (\mbox{$\rho_{t}(s):=\exp(st)$}) distances, Mexican hat wavelets~\citep{HOU2012} (\mbox{$\rho(s):=s^{1/2}\exp(s^{2})$}), and random walks~\citep{SINHA2013}. The behaviour of the filter~$\rho$, and in particular the convergence of \mbox{$1/\rho$} to zero, determines the main properties of the resulting spectral distances, such as smoothness, encoding of local and global shape properties (e.g., Gaussian curvature, geodesic distance), localisation in time and frequency~\citep{HAMMOND2010}. Through a proper filter design and a learning process~\citep{AFLALO2011,BOSCAINI2015A}, we can obtain spectral kernels and distances that are also invariant with respect to isometric transformations or discriminative of a given shape class.

Recalling the definition of the spectral operator in~\citep{PATANE-CGF2017,PATANE2014-PRL,PATANE-STAR2016}, we review equivalent representations of the spectral kernels and distances in terms of the spectral norm of the~$\delta$-functions, the Laplacian spectrum, the spectral operator, and the spectral kernel (Sect.~\ref{sec:CONTINUOUS-SPECTRAL-DIST}). Through these relations, the Laplacian spectral kernels and distances can be computed by approximating the filter with a rational function and converting their evaluation to the solution of a set of differential equations that involve only the Laplace-Beltrami operator. 

The proposed discrete approximation of the spectral kernels (Sect.~\ref{sec:DISCRETE-SPECTRAL-DIST}) is equivalent to solving~$r$ sparse, symmetric, and well-conditioned linear systems, where~$r$ is the degree of the rational polynomial approximation of the filter. Applying an iterative linear solver, the computation of the spectral kernel takes \mbox{$\mathcal{O}(r\tau(n))$} time, where \mbox{$\tau(n)$} typically varies from \mbox{$\tau(n)=n$} to \mbox{$\tau(n)=n\log n$} and depends on the number~$n$ of shape samples and on the sparsity of the coefficient matrix. In a similar way, the computation of the spectral distance between two points has the same order of complexity of the evaluation of the corresponding spectral kernel at the same input points.

Through the spectrum-free approximation, the computed spectral kernel is not affected by the Gibbs phenomenon, as a consequence of the high accuracy of the rational approximation with respect to the low-pass filter associated with the truncated spectral approximation. Furthermore, the spectrum-free approximation is free of user-defined parameters (e.g., the number of Laplacian eigenpairs).

Finally (Sects.~\ref{sec:COMPARISON},~\ref{sec:CONCLUSIONS}), we experimentally verify the accuracy, efficiency, and numerical robustness of the definition and spectrum-free computation of the spectral kernels and distances.

\subsection*{Novelties with respect to previous work}
With respect to previous work and our recent results on the definition and computation of discrete spectral distances~\citep{PATANE-CGF2017,PATANE-STAR2016}, the main novelties of this paper are
\begin{itemize}
\item the discussion of the \emph{optimality of the Laplacian spectrum} for the approximation of the spectral kernel (Sect.~\ref{sec:HLH-BASIS});
\item the study of the \emph{relation between the spectral and Green kernels} (Sect.~\ref{sec:RELATION-KERNEL-DISTANCE}), associated with differential operators (Sect.~\ref{sec:SPEC-OP-DIST});
\item the \emph{analysis of the stability of the spectral distances} with respect to the evaluation of the Laplacian spectrum and to multiple Laplacian eigenvalues (Sect.~\ref{sec:LAPLACIAN-STABILITY}). In fact, the evaluation of the Laplacian spectrum is generally affected by small perturbations of the Laplace-Beltrami operator;
\item a \emph{novel spectrum-free computation} of the spectral distances (Sect.~\ref{sec:CONTINUOUS-SPECTRAL-DISTANCES}), which is achieved by properly approximating the associated filter with a rational polynomial and is expressed in terms of the canonical basis;
\item the \emph{computation of the discrete spectral distances} (Sect.~\ref{sec:UNIFIED-COMPUTATION}), which reduces to the solution of a set of sparse, symmetric, and well-conditioned linear systems. This analogy between the continuous and discrete cases shows the generality of the proposed approach and it has not been addressed by previous work. Furthermore, these results are complementary and more general than the spectrum-free approaches presented in~\citep{PATANE-CGF2017,PATANE-STAR2016};
\item a new \emph{upper bound to the conditioning number} of the linear systems associated with the spectrum-free approximation, which confirms the numerical stability and well-conditioned computation of the class of the spectrum-free approaches (Sect.~\ref{sec:DISCRETE-OPER-DISTANCE});
\item new experiments (Sect.~\ref{sec:COMPARISON}), which enrich the tests initially presented in~\citep{PATANE-CGF2017,PATANE-STAR2016} and address the \emph{robustness} of the computation of the Laplacian spectral distances and kernels with respect to the data resolution, partial sampling, geometric or topological noise, and deformations.
\end{itemize}
\section{Laplacian spectral distances\label{sec:CONTINUOUS-SPECTRAL-DIST}}
Firstly, we present the Laplace-Beltrami operator, its spectrum, and the optimality of the Laplacian eigenbasis (Sect.~\ref{sec:HLH-BASIS}). Recalling the definition of the spectral operator and kernel introduced in~\citep{PATANE-CGF2017,PATANE2014-PRL,PATANE-STAR2016}, we review equivalent representations of the spectral distances in terms of the spectral norm of the~$\delta$-functions, the Laplacian spectrum, the spectral operator, and the spectral kernel (Sect.~\ref{sec:SPEC-OP-DIST}). As novel contribution, we discuss the stability of the spectral distances with respect to the evaluation of the Laplacian spectrum and to multiple Laplacian eigenvalues, which are generally affected by small perturbations of the Laplace-Beltrami operator (Sect.~\ref{sec:LAPLACIAN-STABILITY}). Finally, we introduce a novel relation between the spectral and Green kernels (Sect.~\ref{sec:RELATION-KERNEL-DISTANCE}) and propose a novel spectrum-free computation of the spectral distances in the continuous case, which is presented in a way analogous to the discrete case~\citep{PATANE-STAR2016,PATANE-CGF2017} and shows the generality of the proposed approach (Sect.~\ref{sec:CONTINUOUS-SPECTRAL-DISTANCES}).

\subsection{Laplace-Beltrami operator\label{sec:HLH-BASIS}}
Let~$\mathcal{N}$ be a smooth surface, possibly with boundary, equipped with a Riemannian metrics and let us consider the \emph{inner product} \mbox{$\langle f,g\rangle_{2}:=\int_{\mathcal{N}}fg\textrm{d}\mu$} defined on the space \mbox{$\mathcal{L}^{2}(\mathcal{N})$} of square integrable functions on~$\mathcal{N}$ and the corresponding norm \mbox{$\|\cdot\|_{2}$}. The \emph{Laplace-Beltrami operator}~$\Delta$ is \emph{self-adjoint} \mbox{$\langle\Delta f,g\rangle_{2}=\langle f,\Delta g\rangle_{2}$}, and \emph{positive semi-definite} \mbox{$\langle\Delta f,f\rangle_{2}\geq 0$}, \mbox{$\forall f,g$},~\citep{ROSENBERG1997}.

\paragraph*{Laplacian eigenbasis and its optimality\label{sec:LAPL-EIGIS}}
Recalling that the Laplace-Beltrami operator is self-adjoint and positive semi-definite, it has an \emph{orthonormal eigensystem} \mbox{$(\lambda_{n},\phi_{n})_{n=0}^{+\infty}$}, \mbox{$\Delta\phi_{n}=\lambda_{n}\phi_{n}$}, \mbox{$\lambda_{0}=0$}, \mbox{$\lambda_{n}\leq\lambda_{n+1}$}, in \mbox{$\mathcal{L}^{2}(\mathcal{N})$}. Expressing a function~$f$ in \mbox{$\mathcal{L}^{2}(\mathcal{N})$} in terms of the Laplacian eigenfunctions, we have that
\begin{equation*}
\left\{
\begin{array}{ll}
f=\sum_{n=0}^{+\infty}\alpha_{n}\phi_{n},
&\alpha_{n}:=\langle f,\phi_{n}\rangle_{2};\\
\|f\|_{2}^{2}=\sum_{n=0}^{+\infty}\alpha_{n}^{2},
&\|\nabla f\|_{2}^{2}=\sum_{n=0}^{+\infty}\alpha_{n}^{2}\lambda_{n}.
\end{array}
\right.
\end{equation*}
According to~\citep{AFLALO2015}, the Laplacian eigenfunctions are an optimal basis for the representation of signals with bounded gradient magnitude. In fact, the spectral decomposition \mbox{$f_{n}=\sum_{i=0}^{n}\langle f,\phi_{i}\rangle_{2}\phi_{i}$} is optimal for the approximation of functions with~$\mathcal{L}^{2}$ bounded gradient magnitude; i.e., the residual error \mbox{$r_{n}:=f-f_{n}$} is bounded as
\begin{equation}\label{eq:LS-SPECTRAL-ERROR}
\|r_{n}\|_{2}^{2}
\leq\frac{\|\nabla f\|_{2}^{2}}{\lambda_{n+1}}.
\end{equation}
The spectral decomposition is also optimal for the approximation of functions with respect to the error estimation in (\ref{eq:LS-SPECTRAL-ERROR}). In fact, for any \mbox{$0\leq\alpha<1$} there is no integer~$n$ and no sequence \mbox{$(\psi_{i})_{i=0}^{n}$} of linearly independent functions in~$\mathcal{L}^{2}$ such that
\begin{equation*}\label{eq:LAPL-OPT-BOUND}
\left\|f-\sum_{i=1}^{n}\langle f,\psi_{i}\rangle_{2}\psi_{i}\right\|_{2}
\leq\alpha\frac{\|\nabla f\|_{2}}{\lambda_{n+1}},
\qquad\forall f.
\end{equation*}
Indeed, the Laplacian eigenfunctions are an optimal basis for the definition of scalar functions on a given domain, such as the Laplacian spectral kernels and distances on 3D shapes (Sect.~\ref{sec:DISCRETE-SPECTRAL-DIST}). From the computational point of view, they have two main drawbacks: a high computational cost and storage overhead, which prevent the evaluation of a large number of Laplacian eigenpairs, and numerical instabilities with respect to the surface discretisation~\citep{PATANE-STAR2016,PATANE-CGF2017}.

\subsection{Theoretical background: spectral kernels and distances\label{sec:SPEC-OP-DIST}}
We briefly review the equivalent definitions of the spectral kernels (Sect.~\ref{sec:SPECTRAL-OPERATOR-KERNEL}) and distances (Sect.~\ref{sec:OVERVIEW-SPECTRAL-DISTANCES}), which will be useful to introduce their spectrum-free computation and main properties.

\subsubsection{Spectral operator and kernel\label{sec:SPECTRAL-OPERATOR-KERNEL}}
According to~\citep{PATANE-STAR2016,PATANE-CGF2017}, let \mbox{$\rho:\mathbb{R}^{+}\rightarrow\mathbb{R}$} be a positive \emph{filter map} that is square integrable and admits the power series \mbox{$\rho(s)=\sum_{n=0}^{+\infty}\alpha_{n}s^{n}$}. Considering the orthonormal Laplacian eigensystem \mbox{$(\lambda_{n},\phi_{n})_{n=0}^{+\infty}$}, the spectral representation of the functions \mbox{$\Delta^{i}f$} and \mbox{$(\Delta^{\dag})^{i}f$} is
\begin{equation}\label{eq:SPECTRAL-PSEUDO-INVERSE}
\begin{split}
\Delta^{i}f=\sum_{n=0}^{+\infty}\lambda_{n}^{i}\langle f,\phi_{n}\rangle_{2}\phi_{n},\qquad
(\Delta^{\dag})^{i}f=\sum_{n=1}^{+\infty}\frac{1}{\lambda_{n}^{i}}\langle f,\phi_{n}\rangle_{2}\phi_{n}.
\end{split}
\end{equation}
If \mbox{$\lambda_{k}=0$}, then we neglect the corresponding entry in \mbox{$(\Delta^{\dag})^{i}f$}. In \mbox{$\mathcal{L}^{2}(\mathcal{N})$}, we define the \emph{spectral operator}
\begin{equation}\label{eq:FUNCT-OPER}
\Phi_{\rho}f
=\sum_{n=0}^{+\infty}\alpha_{n}\Delta^{n}f
=\sum_{n=0}^{+\infty}\rho(\lambda_{n})\langle f,\phi_{n}\rangle_{2}\phi_{n},
\end{equation}
which is linear, continuous, self-adjoint, and \mbox{$\Phi_{\rho}f=\langle K_{\rho},f\rangle_{2}$}, with 
\begin{equation}\label{eq:SPECTRAL-KERNEL}
K_{\rho}(\mathbf{p},\mathbf{q})=\sum_{n=0}^{+\infty}\rho(\lambda_{n})\phi_{n}(\mathbf{p})\phi_{n}(\mathbf{q})
\end{equation}
\emph{spectral kernel}. The spectral kernel operator and kernel are well-posed if the filter is bounded or square-integrable. The spectral kernel is also symmetric and self-adjoint, as a consequence of its definition. From these relations, it follows that 
\begin{itemize}
\item the spectral kernel \mbox{$K_{\rho}(\mathbf{p},\cdot)=\Phi_{\rho}\delta_{\mathbf{p}}$} is achieved as the action of the spectral operator on the~$\delta$-function at~$\mathbf{p}$;
\item the spectral kernel is the solution to the differential equation \mbox{$\Phi_{1/\rho}K_{\rho}(\mathbf{p},\cdot)=\delta_{\mathbf{p}}$}. In Sect.~\ref{sec:RELATION-KERNEL-DISTANCE}, this property is further investigated in terms of the relation between the spectral and Green kernels.
\end{itemize}
To express the spectral distances in terms of the spectral operator (Sect.~\ref{sec:CONTINUOUS-SPECTRAL-DISTANCES}), we show that the pseudo-inverse of~$\Phi_{\rho}$ is induced by the filter function~$1/\rho$; i.e., \mbox{$\Phi_{\rho}^{\dag}=\Phi_{1/\rho}$}. In fact, from the spectral representation (\ref{eq:FUNCT-OPER}) of~$\Phi_{1/\rho}$ and~$\Phi_{\rho}$, we get that 
\begin{equation}\label{eq:SPECTRAL-PSEUDOINVERSE}
\begin{array}{l}
\Phi_{\rho}\Phi_{1/\rho}\Phi_{\rho}=\Phi_{\rho},\quad
\Phi_{1/\rho}\Phi_{\rho}\Phi_{1/\rho}=\Phi_{1/\rho},\\
\langle\Phi_{1/\rho}\Phi_{\rho}f,g\rangle_{2}=\langle f,\Phi_{1/\rho}\Phi_{\rho}g\rangle_{2}.
\end{array}
\end{equation}
In the following, we assume that both~$\Phi_{\rho}$ and~$\Phi_{1/\rho}$ are well-defined; for instance, this hypothesis is satisfied if~$\rho$ is not null only on a compact interval of~$\mathbb{R}^{+}$ and it is valid in the discrete case (Sect.~\ref{sec:DISCRETE-SPECTRAL-DIST}), where the Laplacian spectrum belongs to the interval \mbox{$\mathcal{I}:=[0,\lambda_{\max}(\tilde{\mathbf{L}})]$}, with \mbox{$\lambda_{\max}(\tilde{\mathbf{L}})$} maximum Laplacian eigenvalue. 

\subsubsection{Spectral distances\label{sec:OVERVIEW-SPECTRAL-DISTANCES}}
\begin{figure}[t]
\centering
\begin{tabular}{cccc}
\includegraphics[height=70pt]{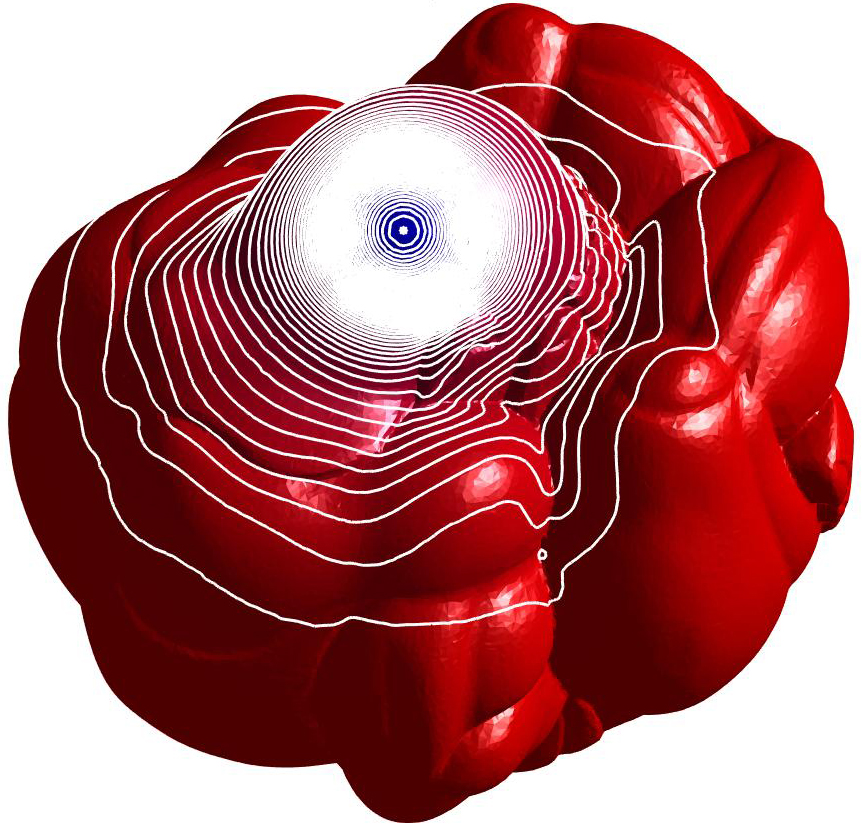}
&\includegraphics[height=70pt]{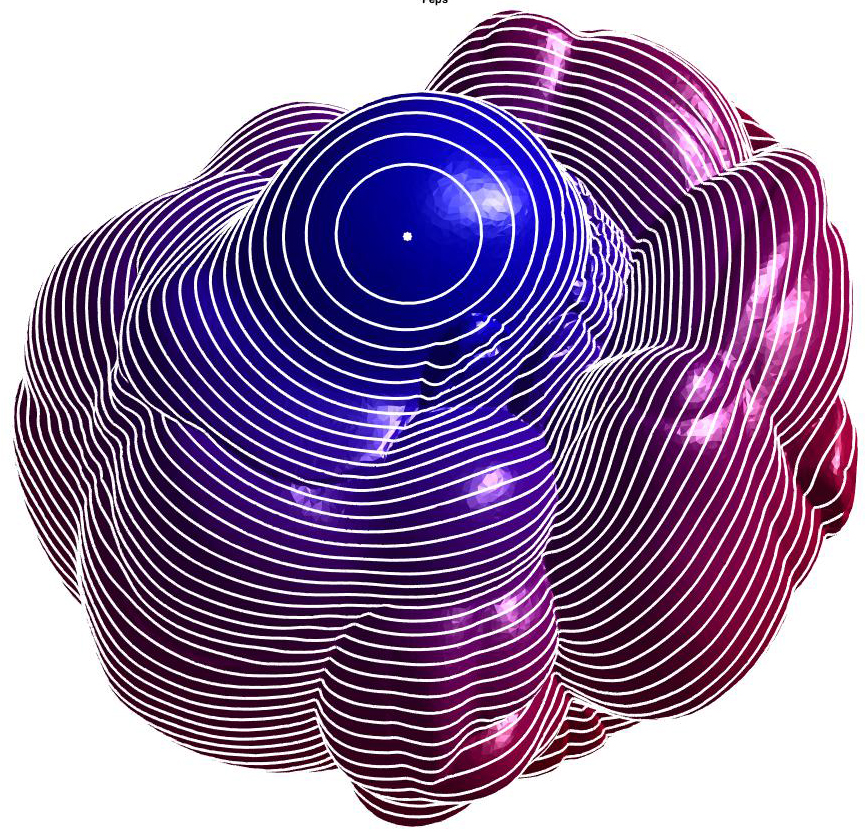}
&\includegraphics[height=70pt]{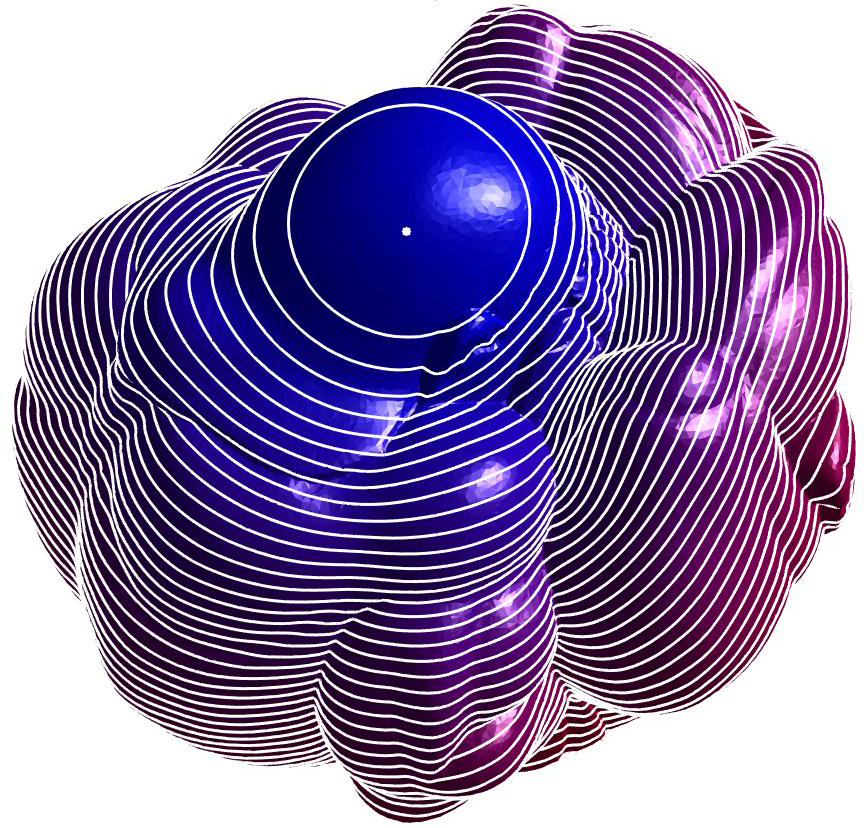}
&\includegraphics[height=70pt]{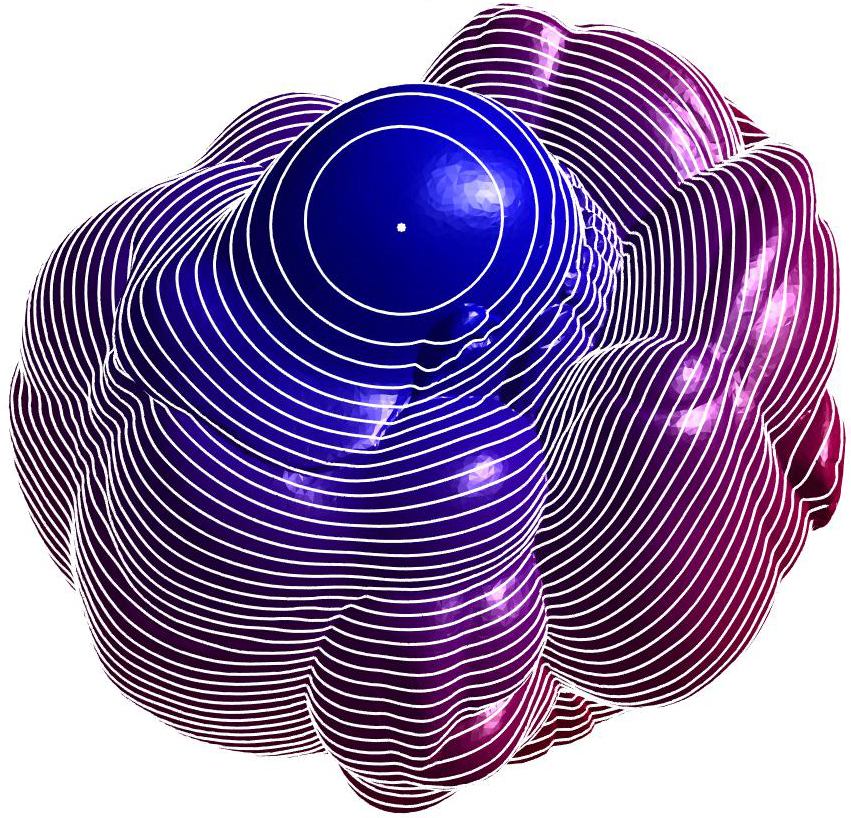}\\
$\rho_{t}(s)=s$ &$\rho_{t}(s)=s^{2}$
&$\rho(s)=s^{3}$ &$\rho(s)=s^{2}\log(1+s)$
\end{tabular}
\caption{Different filters~$\rho$ induce spectral distances from a seed point with a different behaviour, in terms of the locality, shape, and distribution of the level-sets and according to the decay of~$1/\rho$ to zero. These distances have been computed with the Pad\'e-Chebyshev approximation (\mbox{$r=5$}).\label{fig:DIFF-EQUATION}}
\end{figure}

Through the spectral operator, we introduce the scalar product and the corresponding distance as~\citep{PATANE-STAR2016,PATANE-CGF2017} 
\begin{equation}\label{eq:FILT-SCAL-PROD}
\left\{
\begin{array}{ll}
\langle f,g\rangle
:=\langle \Phi_{1/\rho}f,\Phi_{1/\rho}g\rangle_{2}
=\sum_{n=0}^{+\infty}\frac{\langle f,\phi_{n}\rangle_{2}\langle g,\phi_{n}\rangle_{2}}{\rho^{2}(\lambda_{n})}, &(a)\\
d^{2}(f,g)
:=\|f-g\|^{2}
=\sum_{n=0}^{+\infty}\frac{\vert\langle f-g,\phi_{n}\rangle_{2}\vert^{2}}{\rho^{2}(\lambda_{n})}. &(b)
\end{array}
\right.
\end{equation}
Indicating with~$\delta_{\mathbf{p}}$ the map that takes value~$1$ at~$\mathbf{p}$ and~$0$ otherwise, and selecting \mbox{$f:=\delta_{\mathbf{p}}$}, \mbox{$g:=\delta_{\mathbf{q}}$} in Eq. (\ref{eq:FILT-SCAL-PROD}b), the \emph{spectral distance} (Fig.~\ref{fig:DIFF-EQUATION}) on~$\mathcal{N}$ is defined as
\begin{equation*}\label{eq:GEN-COMMUT-DIST}
\begin{split}
d^{2}(\mathbf{p},\mathbf{q})
&:=\|\delta_{\mathbf{p}}-\delta_{\mathbf{q}}\|^{2}\\
&=_{\textrm{Eq. }(\ref{eq:FILT-SCAL-PROD}b)}\sum_{n=0}^{+\infty}\frac{\vert\phi_{n}(\mathbf{p})-\phi_{n}(\mathbf{q})\vert^{2}}{\rho^{2}(\lambda_{n})}\\
&=_{\textrm{Eq. }(\ref{eq:FILT-SCAL-PROD}a)}\|\Phi_{1/\rho}(\delta_{\mathbf{p}})-\Phi_{1/\rho}(\delta_{\mathbf{q}})\|_{2}^{2}\\
&=\|K_{1/\rho}(\mathbf{p},\cdot)-K_{1/\rho}(\mathbf{q},\cdot)\|_{2}^{2}\\
&=K_{1/\rho}(\mathbf{p},\mathbf{p})-2K_{1/\rho}(\mathbf{p},\mathbf{q})+K_{1/\rho}(\mathbf{q},\mathbf{q});
\end{split}
\end{equation*}
i.e., these equivalent formulations involve the Laplacian spectrum, the spectral operator and kernel. The third equality follows from the identity \mbox{$\Phi_{\rho}(\delta_{\mathbf{p}})=K_{\rho}(\mathbf{p},\cdot)$} and will be applied to the computation of the spectral distances (Sect.~\ref{sec:UNIFIED-COMPUTATION}); in fact, it is independent of the evaluation of the Laplacian spectrum. The last identity is achieved by applying the relation
\begin{equation*}\label{eq:SPECTRAL-DISTANCE-SQUARE}
\begin{split}
d^{2}(\mathbf{p},\mathbf{q})
&=\|K_{1/\rho}(\mathbf{p},\cdot)-K_{1/\rho}(\mathbf{q},\cdot)\|_{2}^{2}\\
&=\sum_{n=0}^{+\infty}\rho^{-2}(\lambda_{n})\phi_{n}(\mathbf{p})\phi_{n}(\mathbf{p})
-2\sum_{n=0}^{+\infty}\rho^{-2}(\lambda_{n})\phi_{n}(\mathbf{p})\phi_{n}(\mathbf{q})\\
&+\sum_{n=0}^{+\infty}\rho^{-2}(\lambda_{n})\phi_{n}(\mathbf{q})\phi_{n}(\mathbf{q}).
\end{split}
\end{equation*}
\subsection{Stability of the eigenpairs of the spectral operator\label{sec:LAPLACIAN-STABILITY}}
Generalising the results in~\citep{PATANE-STAR2016,PATANE-CGF2017}, we show that the computation of a single eigenvalue of the spectral operator is numerically stable and instabilities are generally due to repeated or close eigenvalues. Firstly, we notice that if~$\lambda$ is a Laplacian eigenvalue of multiplicity~$m$ then \mbox{$\mu:=\rho(\lambda)$} is an eigenvalue of~$\Phi_{\rho}$ and its multiplicity is equal to or greater than~$m$. The corresponding eigenfunction is the Laplacian eigenfunction associated with the eigenvalue~$\lambda$; i.e., \mbox{$\Delta\phi=\lambda\phi$} and \mbox{$\Phi_{\rho}\phi=\rho(\lambda)\phi$}. 

We perturb the spectral operator by~$\delta \mathcal{E}$, \mbox{$\delta\rightarrow 0$}, and compute the eigenpair \mbox{$(\mu(\delta),\phi(\delta))$} of the corresponding operator \mbox{$\Phi_{\rho}+\delta\mathcal{E}$}; i.e.,
\begin{equation}\label{eq:PERTURBED-SPECT-OP}
\left(\Phi_{\rho}+\delta\mathcal{E}\right)\phi(\delta)=\mu(\delta)\phi(\delta),\qquad
\phi(0)=\phi, \quad \mu(0)=\mu.
\end{equation}
Recalling that the derivative of \mbox{$\mu(\delta)$} measures the variation of the eigenvalue, deriving (\ref{eq:PERTURBED-SPECT-OP}) with respect to~$\delta$, and evaluating the resulting relation at~$0$, we have that
\begin{equation}\label{eq:SPECTRAL-RELATION}
\mathcal{E}\phi+\Phi_{\rho}\phi^{\prime}(0)=\mu^{\prime}(0)\phi+\mu\phi^{\prime}(0).
\end{equation}
From (\ref{eq:SPECTRAL-RELATION}) and the self-adjointness of~$\Phi_{\rho}$, it follows that \mbox{$\langle\phi,\Phi_{\rho}\phi^{\prime}(0)\rangle_{2}=\mu\langle\phi,\phi^{\prime}(0)\rangle_{2}$}. Multiplying both sides of (\ref{eq:SPECTRAL-RELATION}) by~$\phi$ and applying the previous identity we get that
\begin{equation*}
\vert\mu^{\prime}(0)\vert
=\vert\langle\phi,\mathcal{E}\phi\rangle_{2}\vert
\leq\|\mathcal{E}\|_{2}\|\phi\|_{2}^{2}=\|\mathcal{E}\|_{2};
\end{equation*}
i.e., the computation of the eigenvalue of~$\Phi_{\rho}$ with multiplicity one is stable. 

Assuming that \mbox{$\mu:=\rho(\lambda)$} is an eigenvalue of~$\Phi_{\rho}$ with multiplicity~$m$ and rewriting the characteristic polynomial as \mbox{$p(s)=(s-\mu)^{m}q(s)$}, where \mbox{$q(\cdot)$} is a polynomial of degree \mbox{$n-m$} and \mbox{$q(\mu)\neq 0$}, we get that 
\begin{equation*}
(s-\mu)^{m}=\frac{p(s)}{q(s)}
\approx \frac{O(\delta)}{q(s)},\quad
\delta\rightarrow 0,\quad
\textrm{as }
p(s)\rightarrow 0,
\textrm{ when } s\rightarrow\mu;
\end{equation*}
i.e., \mbox{$s\approx\mu+O(\delta^{\frac{1}{m}})$}. Indeed, modifying the Laplacian matrix in such a way that the filtered eigenvalues are perturbed by \mbox{$\delta:=10^{-m}$} corresponds to a change of order~$0.1$ in \mbox{$\mu$} (i.e., \mbox{$s\approx\mu+0.1$}) and this amplification becomes larger as the multiplicity of the eigenvalue increases.

While multiple eigenvalues are typically associated with symmetric shapes, numerically close or switched eigenvalues are present regardless of the surface regularity. This unstable computation of multiple eigenpairs of the spectral operator generally affects the accuracy of the \emph{truncated spectral approximation}
\begin{equation}\label{eq:TRUNCATED-SPECTRAL-APPROXIMATION}
d_{k}(\mathbf{p},\mathbf{q})
:=\sum_{n=0}^{k}\frac{\vert\phi_{n}(\mathbf{p})-\phi_{n}(\mathbf{q})\vert^{2}}{\rho^{2}(\lambda_{n})}
\end{equation}
of the corresponding distances. In fact, each filtered eigenvalue \mbox{$\rho(\lambda_{n})$}, \mbox{$n=1,\ldots,k$}, which appears at the denominator of~$d_{k}$, can further accentuate the numerical error of~$\lambda_{n}$ in the distance computation. Furthermore, the computation of the Laplacian spectrum is time-consuming and it is difficult to properly select the number of eigenpairs that is necessary to accurately approximate the spectral distance. Indeed, we propose an evaluation of the spectral distance that is independent of the computation of the Laplacian spectrum and is equivalent to a set of differential equations involving only the Laplace-Beltrami operator and its pseudo-inverse. This novel spectrum-free computation is a generalisation of the discrete approach recently presented in~\citep{PATANE-STAR2016,PATANE-CGF2017}.
\subsection{Relation between spectral and Green kernels\label{sec:RELATION-KERNEL-DISTANCE}}
To study the relation between the Green and the spectral kernels, let us consider a linear differential operator~$\mathcal{L}$ and the corresponding \emph{Green kernel} \mbox{$K:\mathcal{N}\times\mathcal{N}\rightarrow\mathbb{R}$} such that \mbox{$\mathcal{L}K(\mathbf{p},\mathbf{q})=\delta(\mathbf{p}-\mathbf{q})$}. Then, the solution to the differential equation \mbox{$\mathcal{L}u=f$} is expressed in terms of the Green kernel as \mbox{$u(\mathbf{p})=\langle K(\mathbf{p},\cdot),f\rangle_{2}$}. 

\begin{figure*}[t]
\centering
\begin{tabular}{ccccccc}
(a)\includegraphics[height=100pt]{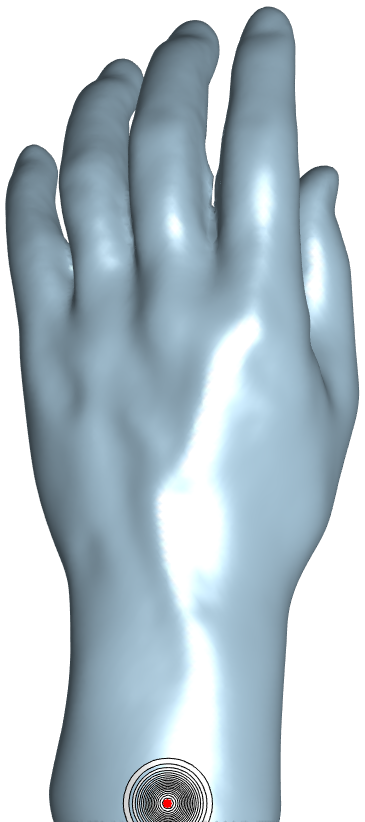}
&\includegraphics[height=100pt]{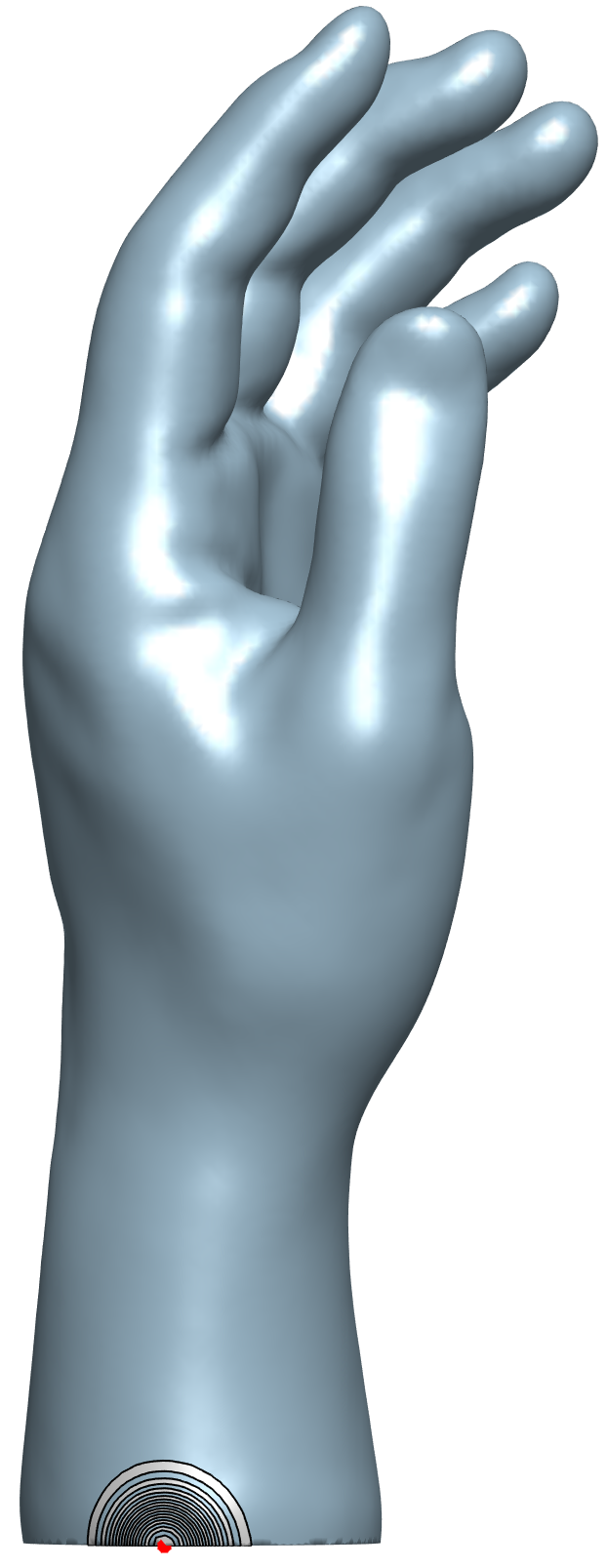}
&\includegraphics[height=100pt]{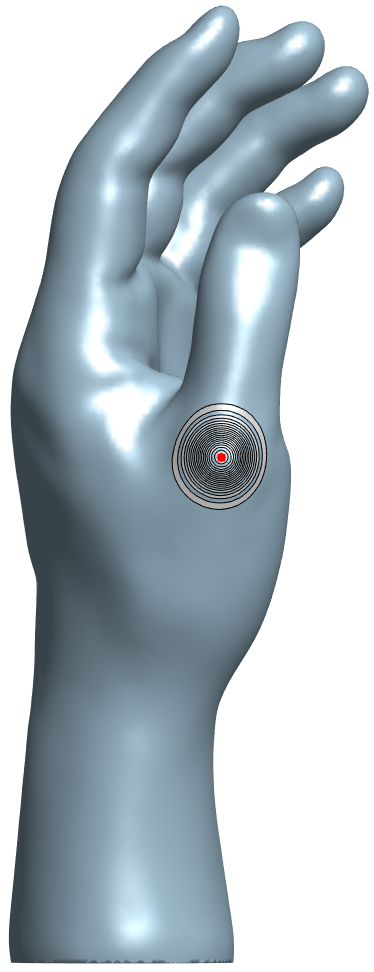}
&\includegraphics[height=100pt]{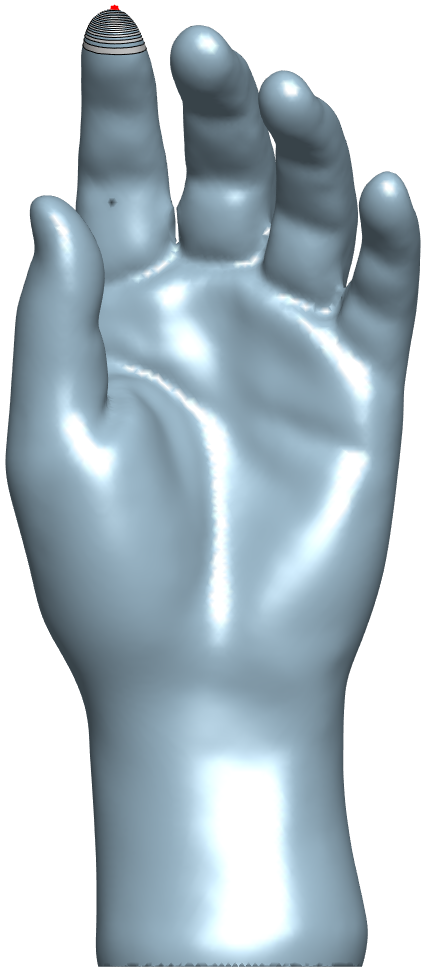}
&\includegraphics[height=100pt]{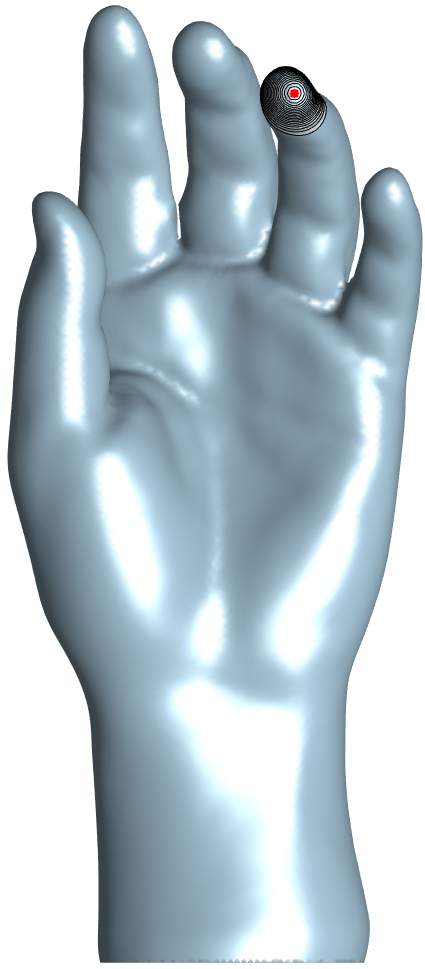}
&\includegraphics[height=100pt]{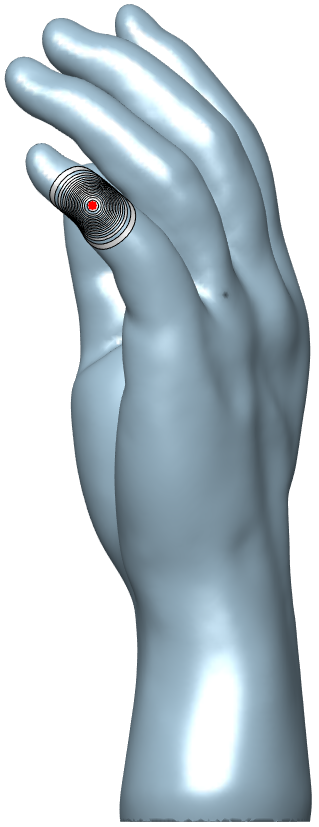}\\
\hline
(b)\includegraphics[height=100pt]{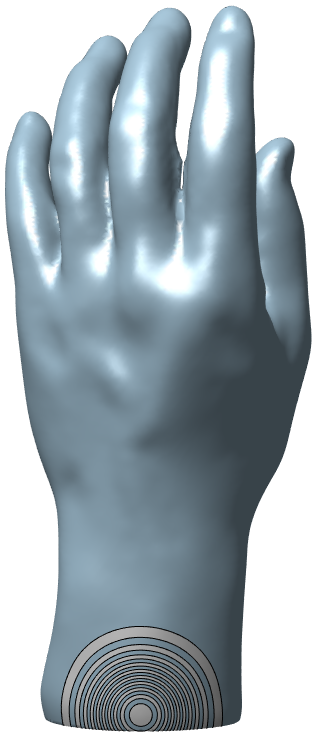}
&\includegraphics[height=100pt]{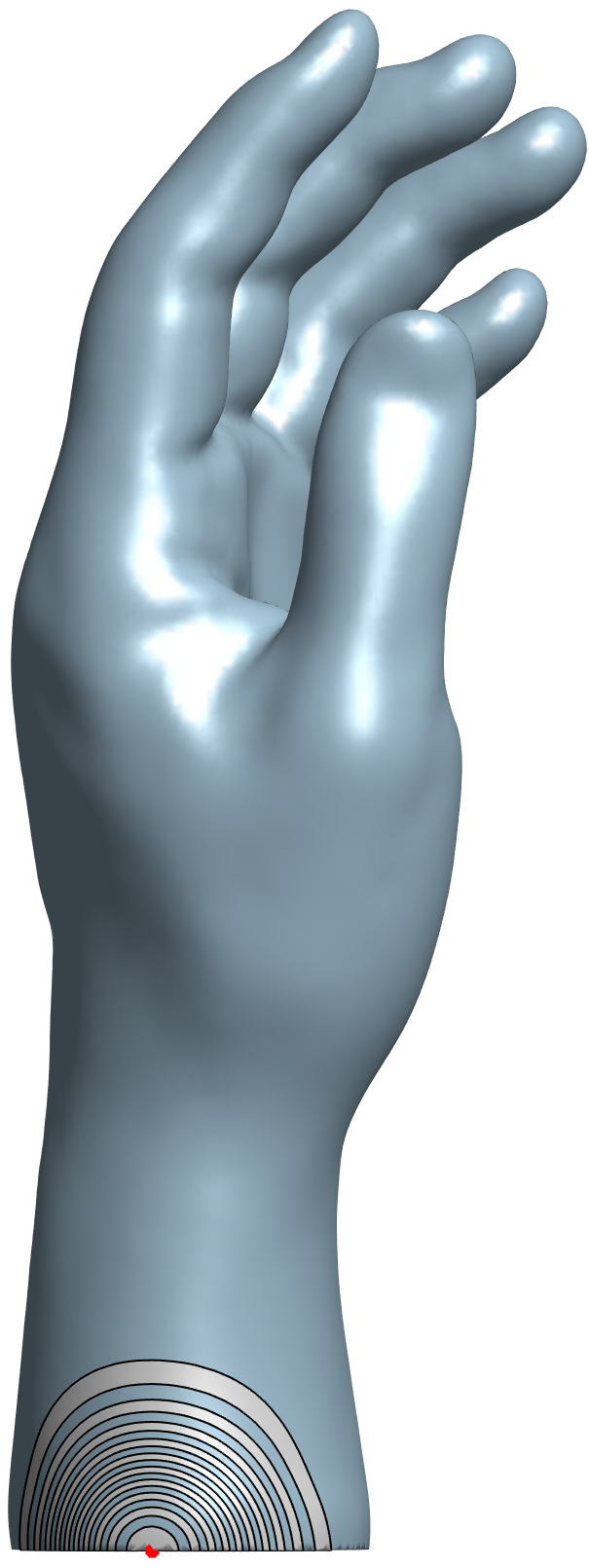}
&\includegraphics[height=100pt]{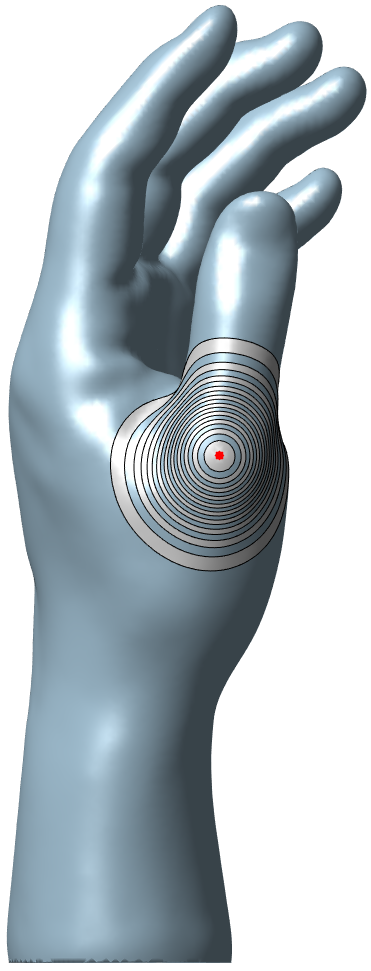}
&\includegraphics[height=100pt]{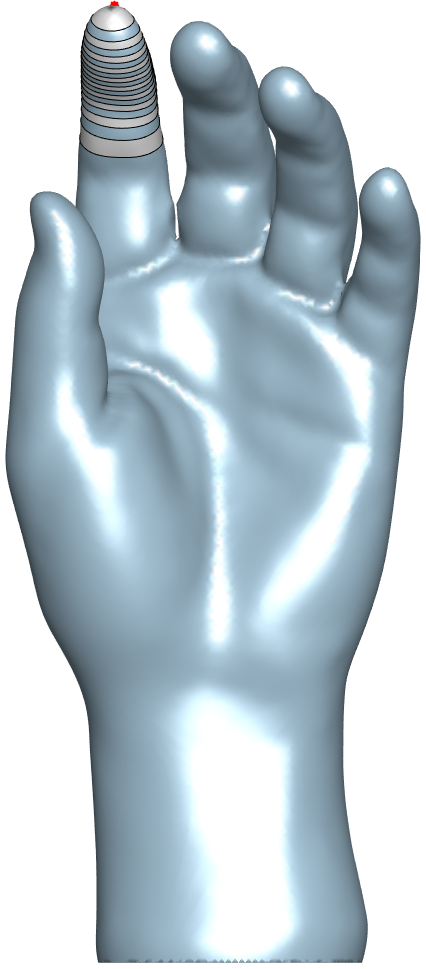}
&\includegraphics[height=100pt]{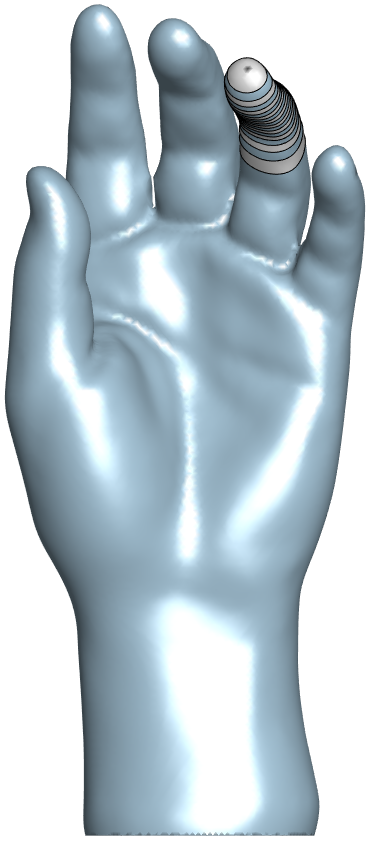}
&\includegraphics[height=100pt]{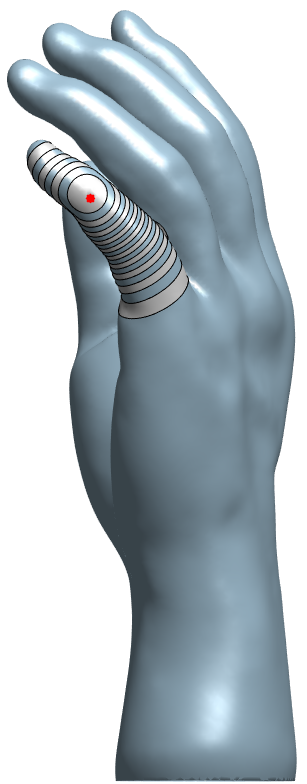}
\end{tabular}
\caption{(a,b) Level-sets of the diffusion kernel at two different scales and at different seed points, which confirm the locality, smoothness, and shape-awareness of the heat kernel.\label{fig:HEAT-DIFFUSION}}
\end{figure*}

Noting that the eigensystem \mbox{$(\mu_{n},\psi_{n})_{n=0}^{+\infty}$} of the integral operator \mbox{$\mathcal{A}_{K}f:=\langle K(\cdot,\cdot),f\rangle_{2}$} induced by the Green kernel satisfies the relation \mbox{$\mathcal{L}\psi_{n}=\mu_{n}^{-1}\psi_{n}$}, \mbox{$\mu_{n}\neq 0$}, we have that the eigensystem of~$\mathcal{L}$ is \mbox{$(\mu_{n}^{-1},\psi_{n})_{n=1}^{+\infty}$}. Indeed,~$\mathcal{L}$ and~$\mathcal{A}_{K}$ have the same eigenfunctions and reciprocal eigenvalues. In particular, the spectral representation of the Green kernel \mbox{$K(\mathbf{p},\mathbf{q})=\sum_{n=0}^{+\infty}\mu_{n}\psi_{n}(\mathbf{p})\psi_{n}(\mathbf{q})$} is uniquely defined by the spectrum of~$\mathcal{L}$. Combining the previous results with the properties of the spectral operator and kernel (Sect.~\ref{sec:OVERVIEW-SPECTRAL-DISTANCES}), we get that
\begin{itemize}
\item for the \emph{harmonic operator} \mbox{$\mathcal{L}=\Delta$}, the Green kernel is the \emph{commute-time kernel}
\begin{equation*}
K_{\Delta}(\mathbf{p},\mathbf{q})=\sum_{n=1}^{+\infty}\frac{1}{\lambda_{n}}\phi_{n}(\mathbf{p})\phi_{n}(\mathbf{q});
\end{equation*}
\item for the \emph{bi-Laplacian operator} \mbox{$\mathcal{L}:=\Delta^{2}$}, the Green kernel is the \emph{bi-harmonic kernel}
\begin{equation*}
K_{\Delta^{2}}(\mathbf{p},\mathbf{q}):=\sum_{n=1}^{+\infty}\frac{1}{\lambda_{n}^{2}}\phi_{n}(\mathbf{p})\phi_{n}(\mathbf{q});
\end{equation*}
\item for the \emph{diffusion operator} \mbox{$\mathcal{L}=\exp(t\Delta)$}, the Green kernel is the \emph{diffusion kernel} (Fig.~\ref{fig:HEAT-DIFFUSION})
\begin{equation*}
K_{t}(\mathbf{p},\mathbf{q})=\sum_{n=0}^{+\infty}\exp(-t\lambda_{n})\phi_{n}(\mathbf{p})\phi_{n}(\mathbf{q});
\end{equation*}
\item for the \emph{spectral operator} \mbox{$\mathcal{L}=\Phi_{1/\rho}$}, the Green kernel is the \emph{spectral kernel}~$K_{\rho}$ defined in Eq. (\ref{eq:SPECTRAL-KERNEL}).
\end{itemize}
\subsection{Spectrum-free approximation of kernels and distances\label{sec:CONTINUOUS-SPECTRAL-DISTANCES}}
Recalling the relation \mbox{$\Phi_{\rho}^{\dag}=\Phi_{1/\rho}$} in Eq. (\ref{eq:SPECTRAL-PSEUDOINVERSE}) and noting that
\begin{equation*}\label{eq:SPECTRAL-DIST-OPERATOR}
d(f,g)=\|u\|_{2},\quad 
u=\Phi_{1/\rho}(f-g) \Longleftrightarrow \Phi_{\rho}u=f-g,
\end{equation*}
the spectral distance is equivalent to the norm of (i) the function \mbox{$u=\Phi_{1/\rho}(f-g)$} and (ii) the solution of the differential equation \mbox{$\Phi_{\rho}u=f-g$}. Through these relations, the spectral kernels will be computed by approximating the filter with a rational function (Sect.~\ref{sec:KERNEL-APPROX-CONV-ACCURACY}) and by converting the evaluation of the corresponding distances (Sect.~\ref{sec:RATIONAL-APPROX}) to the solution of a set of differential equations that involve the Laplace-Beltrami operator. The choice of one of these two equivalent representations depends on the selected filter and the simplicity of evaluating either~$\Phi_{\rho}$ or~$\Phi_{1/\rho}$.
\subsubsection{Kernel approximation, convergence, and accuracy\label{sec:KERNEL-APPROX-CONV-ACCURACY}}
Let us introduce the space of the Laplacian spectral kernels as
\begin{equation*}
\mathcal{K}(\mathcal{N}):=\{K_{\rho}:\mathcal{N}\times\mathcal{N}\rightarrow\mathbb{R},\,K_{\rho} \textrm{ spectral kernel in Eq. (\ref{eq:SPECTRAL-KERNEL})}\},
\end{equation*}
and equipped with the \mbox{$\mathcal{L}^{2}(\mathcal{N}\times\mathcal{N})$} scalar product. Since
\begin{equation*}\label{eq:KERNEL-DISTANCE}
\|K_{\rho}-K_{\varphi}\|_{2}^{2}
=\sum_{n=0}^{+\infty}\vert\rho(\lambda_{n})-\varphi(\lambda_{n})\vert^{2}\\
\leq\|\rho-\varphi\|_{2}^{2},
\end{equation*}
the approximation of a given spectral kernel~$K_{\rho}$ with a new kernel~$K_{\varphi}$ in \mbox{$\mathcal{K}(\mathcal{N})$} is reduced to the approximation of~$\rho$ by~$\varphi$ on a proper subspace of functions (e.g., the space of polynomials or rational polynomials). The class of functions used for the approximation of the input filter~$\rho$ is selected in such a way that~$K_{\varphi}$ provides a good approximation of the input filter~$K_{\rho}$ and is easily computable. 
\begin{figure}[t]
\centering
\begin{tabular}{ccc}
\includegraphics[height=110pt]{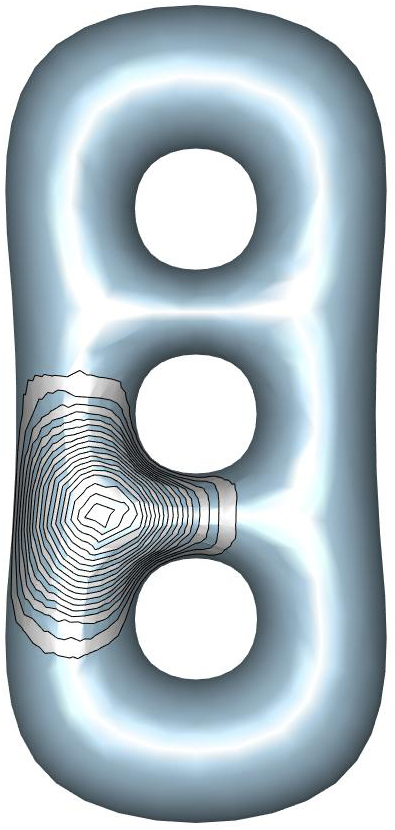}
&\includegraphics[height=110pt]{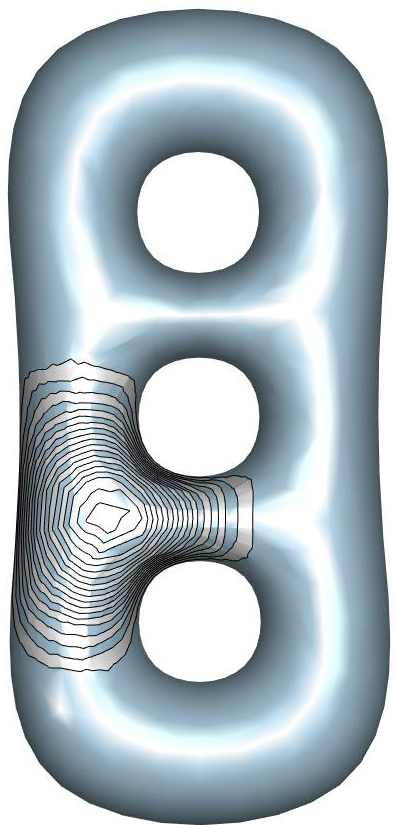}
&\includegraphics[height=110pt]{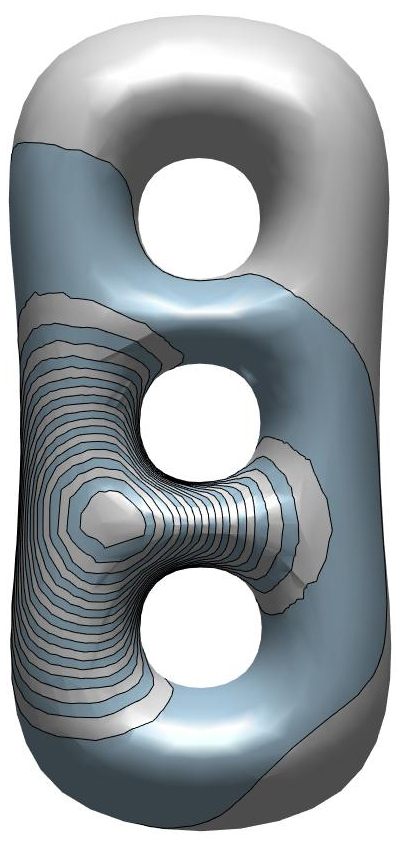}\\
$\tilde{\mathbf{L}}^{2}\mathbf{e}_{i}$ &$\tilde{\mathbf{L}}^{3}\mathbf{e}_{i}$ &$\tilde{\mathbf{L}}^{10}\mathbf{e}_{i}$\\
\includegraphics[height=110pt]{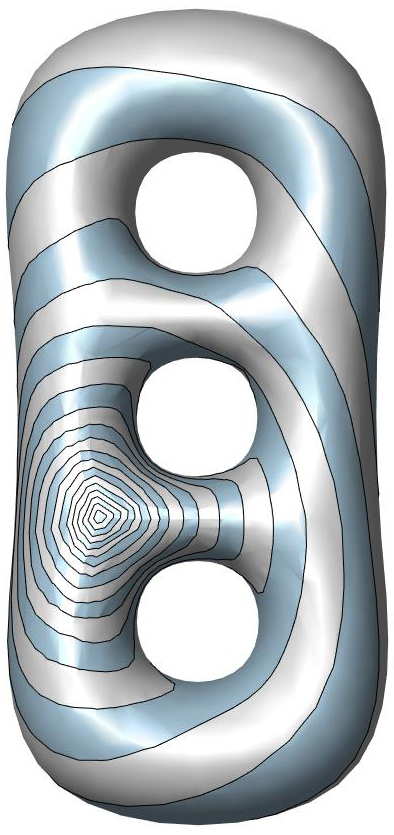}
&\includegraphics[height=110pt]{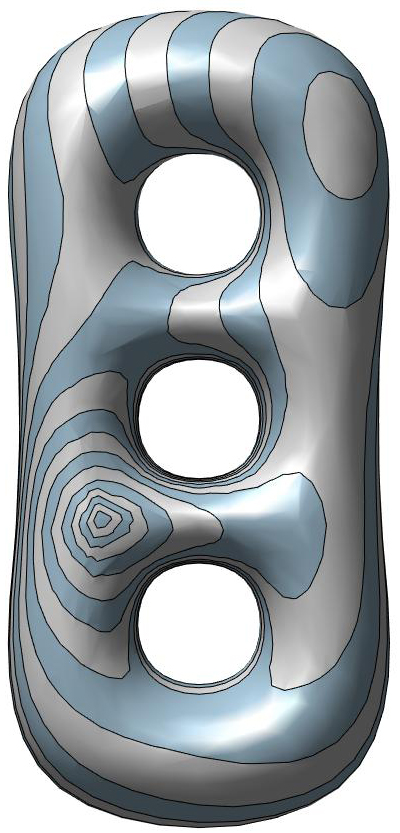}
&\includegraphics[height=110pt]{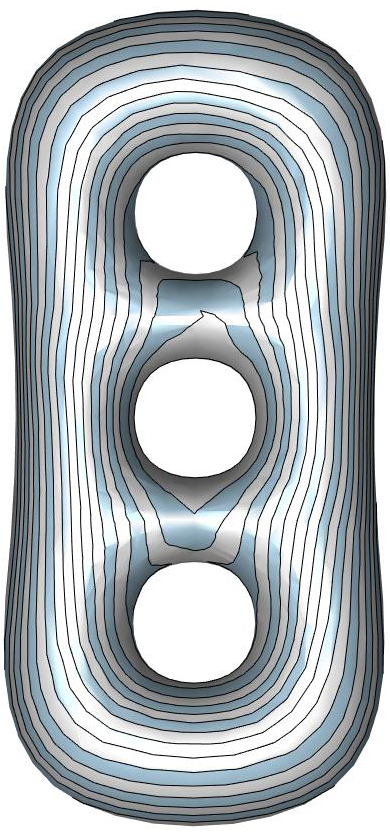}\\
$(\tilde{\mathbf{L}}^{\dag})^{2}\mathbf{e}_{i}$ &($\tilde{\mathbf{L}}^{\dag})^{3}\mathbf{e}_{i}$ &($\tilde{\mathbf{L}}^{\dag})^{10}\mathbf{e}_{i}$\\
\end{tabular}
\caption{Level-sets of the basis functions at~$\mathbf{p}_{i}$ associated with the rational polynomial approximation of the spectral distances.\label{fig:BASIS-FUNCTIONS}}
\end{figure}
\subsubsection{Rational approximation based on the canonical basis\label{sec:RATIONAL-APPROX}}
We apply the \emph{Pad\'e-Chebyshev rational approximation} to the filter map~$1/\rho$. Let~$\mathcal{R}_{r}^{l}$ be the space of all rational functions
\begin{equation*}
c_{rl}(s)
:=\frac{p_{l}(s)}{q_{r}(s)}
=\frac{\beta_{0}+\beta_{1}s+\ldots+\beta_{l}s^{l}}{\alpha_{0}+\alpha_{1}s+\ldots+\alpha_{r}s^{r}},\qquad
s\in[a,b].
\end{equation*}
\begin{figure}[t]
\centering
\begin{tabular}{ccc||c}
\hline
\multicolumn{3}{c||}{\emph{Truncated  spectral approx.}}
&\multicolumn{1}{c}{\emph{P.C. approx.}}\\
\hline
(a)~$k=100$		&(b)~$k=1000$		&(c)~$k=2K$		&(d)~$\epsilon_{\infty}=10^{-6}$\\
\hline\hline
\multicolumn{3}{c||}{$s=10^{-1}$}
&\multicolumn{1}{c}{}\\
\includegraphics[height=90pt]{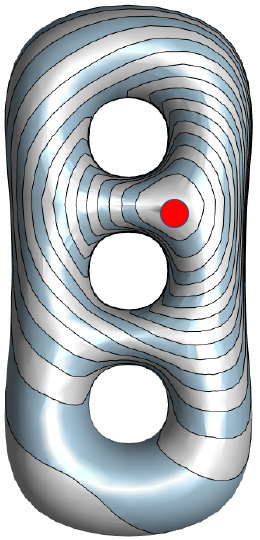}
&\includegraphics[height=90pt]{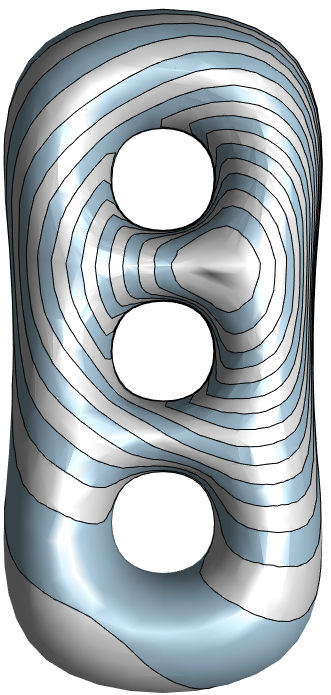}
&\includegraphics[height=90pt]{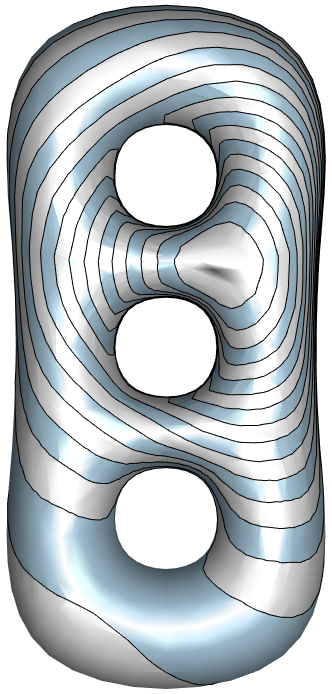}
&\includegraphics[height=90pt]{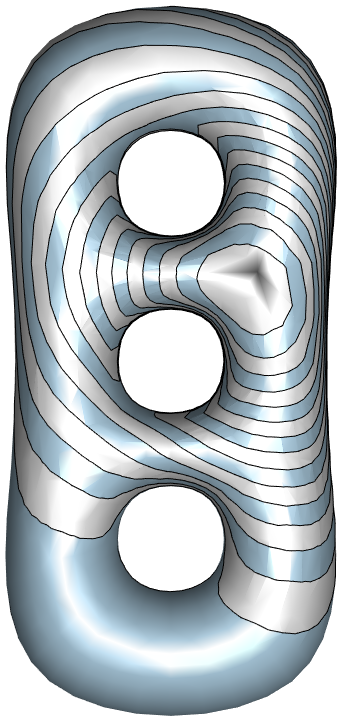}\\
\hline
\multicolumn{3}{c||}{$s=10^{-3}$}
&\multicolumn{1}{c}{}\\
\includegraphics[height=90pt]{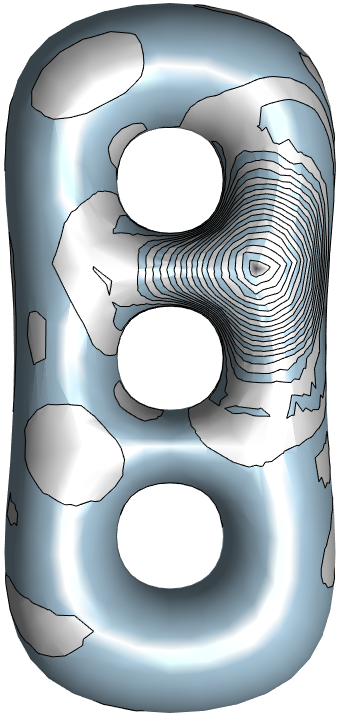}
&\includegraphics[height=90pt]{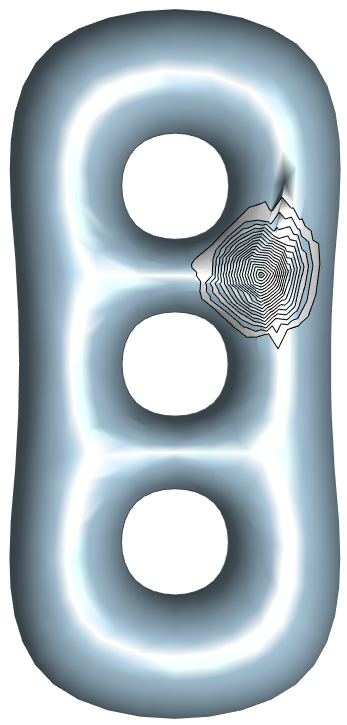}
&\includegraphics[height=90pt]{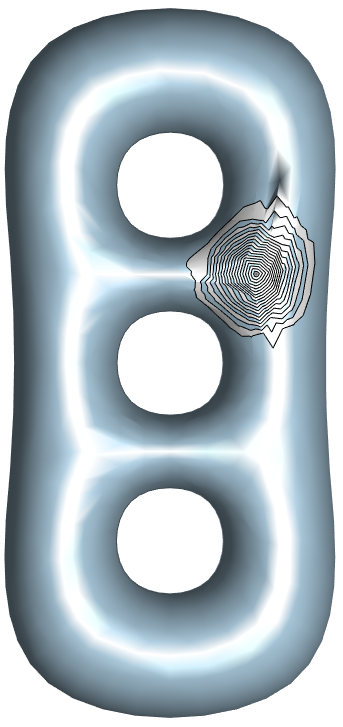}
&\includegraphics[height=90pt]{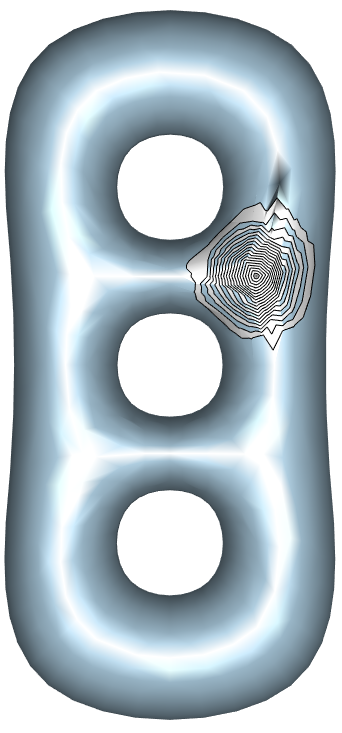}
\end{tabular}
\caption{Smoothness and locality of the diffusion kernel at a seed point (red dot) and at different scales, which have been computed with (a-c) the truncated spectral approximation (i.e.,~$k$ Laplacian eigenpairs) and (d) the Pad\`e-Chebyshev approximation. Since the input shape has~$2K$ vertices, the spectral approximation (c) provides the ground-truth.\label{fig:DIFFUSIVE-COMPUTATION}}
\end{figure}
We briefly recall that the representation of the rational approximation is not unique, unless we impose that it is irreducible; for example, by choosing \mbox{$q_{r}(a)=1$}. Given a filter \mbox{$\rho:\,[a,b]\rightarrow\mathbb{R}$}, there exists a unique best approximation of~$1/\rho$ in~$\mathcal{R}_{r}^{l}$ with respect to the~$\ell_{\infty}$ norm~\citep{GOLUB1989} (Ch.~9), which is represented in terms of the \emph{canonic basis} \mbox{$\mathcal{B}:=\{s^{i}\}_{i=0}^{l}\cup\{s^{-i}\}_{i=1}^{r}$}.
\begin{figure}[t]
\centering
\begin{tabular}{cc}
\includegraphics[height=110pt]{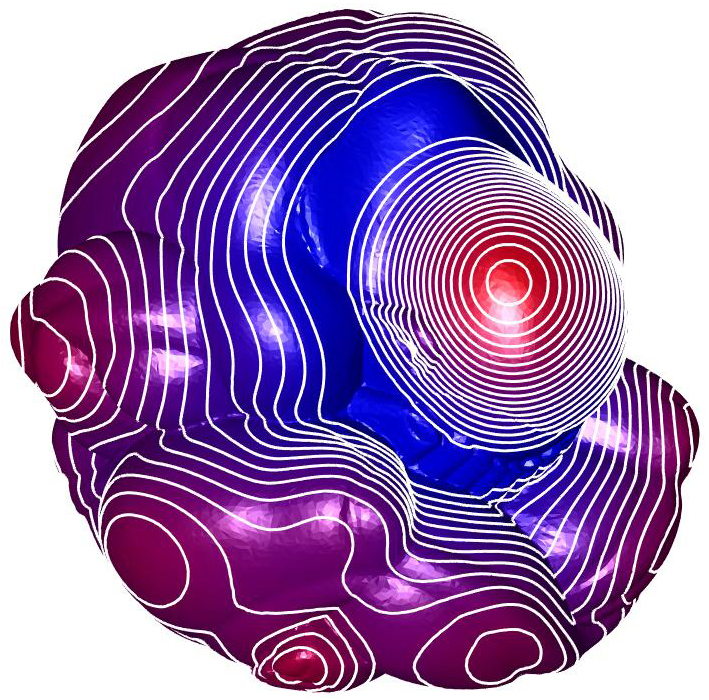}
&\includegraphics[height=110pt]{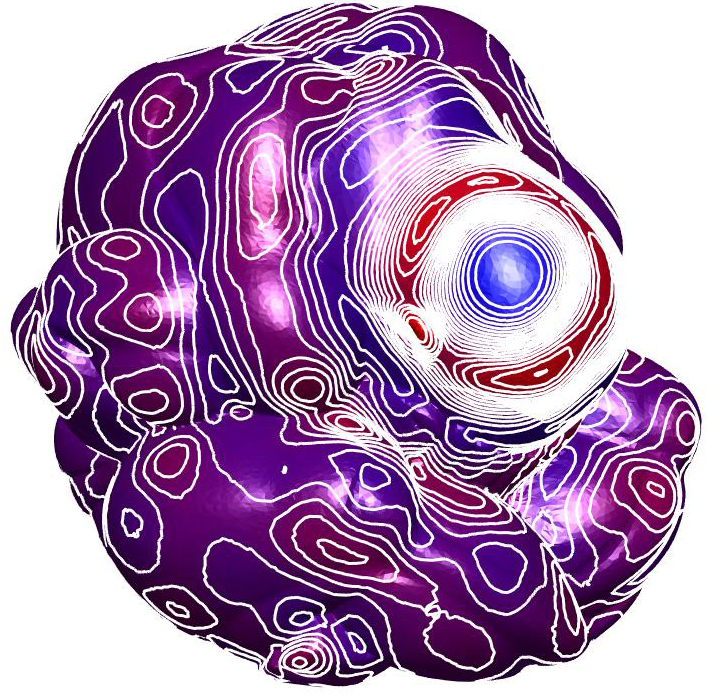}\\
(a)~$\rho(s)=\frac{(1-s)^2}{1+s}$ &(b)~$\rho$, trunc. spect. approx.
\end{tabular}
\caption{Comparison of the spectral distance induced by the same filter and approximated with the (a) spectrum-free and (b) truncated spectral approximation (\mbox{$k=500$}). This last method produces some artefacts that become evident as we move far from the seed point.\label{fig:OMOTONDO-RATIONAL}}
\end{figure}

Recalling that the dimension of the space~$\mathcal{R}_{r}^{l}$ of rational polynomial of degree \mbox{$(r,l)$} is \mbox{$l+r+1$}, let us express the best rational approximation~$c_{rl}$ of~$\rho$ in terms of these basis functions. To this end, we can either apply algebraic rules or impose interpolating constraints at \mbox{$(l+r+1)$} points in order to compute the new coefficients through the identity
\begin{equation*}\label{eq:SPECTRAL-CANONIC-BASIS}
\sum_{i=0}^{l}a_{i}s^{i}+\sum_{i=1}^{r}b_{i}s^{-i}=c_{rl}(s).
\end{equation*}
The canonical basis used to represent the rational approximation of the input filter simplifies the computation of the spectral kernel; in fact, we need to evaluate only the functions \mbox{$\Delta^{i}f$} and \mbox{$(\Delta^{\dag})^{i}f$}. Instead of applying the spectral representation (\ref{eq:SPECTRAL-PSEUDO-INVERSE}), we define a recursive procedure as follows. To compute \mbox{$g_{1}=\Delta^{\dag}f$}, let us multiply both sides by~$\Delta$ and notice that
\begin{equation}\label{eq:CONTINUOUS-APPROX}
\Delta g_{1}
=\Delta\Delta^{\dag}f
=f-\langle f,\phi_{0}\rangle_{2}\phi_{0},
\end{equation}
where~$\phi_{0}=\mathbb{1}$ is the constant eigenfunction equal to~$1$. By definition of pseudo-inverse,~$g_{1}$ is the least-squares solution to the equation \mbox{$\Delta g_{1}=\tilde{f}$}, where~$\tilde{f}$ is equal to~$f$ minus its mean \mbox{$\langle f,\phi_{0}\rangle_{2}\phi_{0}$}. For the general case, we apply the recursive relation
\begin{equation*}
g_{i}
:=(\Delta^{\dag})^{i}f
=\Delta^{\dag}(\Delta^{\dag})^{i-1}f
=\Delta^{\dag}g_{i-1},\quad i\geq 2,
\end{equation*}
which reduces to the previous case. In a similar way, \mbox{$g_{i}=\Delta^{i} f$} is calculated as \mbox{$h_{i}=\Delta h_{i-1}$}, \mbox{$i=1,\ldots,r$}, with \mbox{$h_{0}:=f$}. Indeed, the proposed approach requires the solution of \mbox{$(l+r)$} Laplace equation with a different right-hand side.

\paragraph*{Convergence and accuracy}
To verify that the sequence 
\begin{equation*}
(\Phi_{1/\rho}^{(r)}f)_{r=0}^{+\infty},\quad
\Phi_{1/\rho}^{(r)}f
:=\sum_{n=0}^{+\infty}c_{rl}(\lambda_{n})\langle f,\phi_{n}\rangle_{2}\phi_{n},
\end{equation*}
induced by the rational polynomial approximation~$c_{rl}$ of~$1/\rho$, converges to \mbox{$\Phi_{1/\rho}f$}, we apply the upper bound
\begin{equation*}\label{eq:CONTINUOUS-APPROX-ACCURCACY}
\begin{split}
\left\|\Phi_{1/\rho}^{(r)}f-\Phi_{1/\rho}f\right\|_{2}^{2}
&\leq \|c_{rl}-1/\rho\|_{\infty}^{2}\sum_{n=0}^{+\infty}\vert\langle f,\phi_{n}\rangle_{2}\vert^{2}\\
&=\sigma_{rl}^{2}\|f\|_{2}^{2},\quad\sigma_{rl}\approx\mathcal{O}(s^{r+l+1}),\quad s\rightarrow 0;\\
\end{split}
\end{equation*}
where~$\sigma_{rl}$ is the approximation error between~$1/\rho$ and~$c_{rl}$. For the diffusion operator~\citep{VARGA1990}, \mbox{$\sigma_{rr}=10^{-5}$}.

\section{Spectrum-free approximation of kernels and distances\label{sec:DISCRETE-SPECTRAL-DIST}}
We briefly summarise the discretisation of the spectral kernel and distances~\citep{PATANE-STAR2016,PATANE-CGF2017}, which is then used to establish a simple relation between the spectral kernel and its pseudo-inverse (Sect.~\ref{sec:DISCRETE-OPER-DISTANCE}). Then, we introduce the spectrum-free approximation of the spectral kernels and distances (Sect.~\ref{sec:UNIFIED-COMPUTATION}) and discuss the computational cost of the main steps of the proposed approach (Sect.~\ref{sec:COMPUTATIONAL-COST}).
\begin{figure}[t]
\centering
\begin{tabular}{ccc}
\includegraphics[width=100pt]{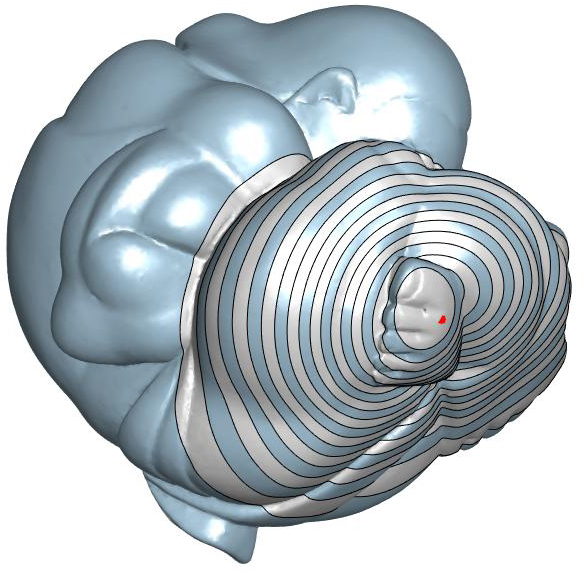}
&\includegraphics[width=100pt]{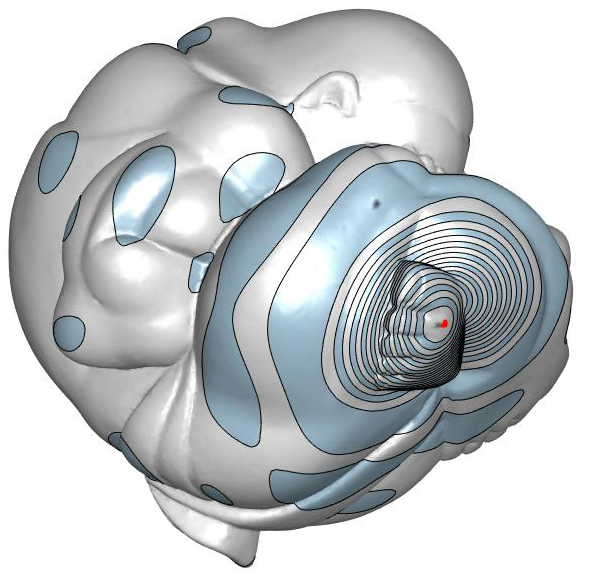}
&\includegraphics[width=100pt]{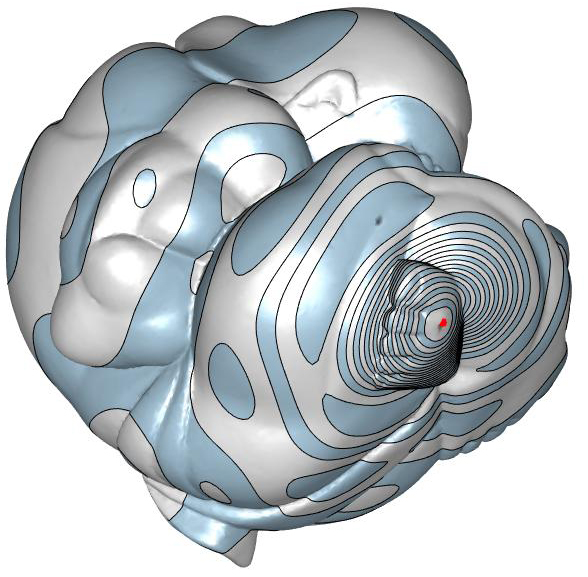}\\
$t=10^{-1}$ &$t=10^{-2}$ &$t=10^{-3}$
\end{tabular}
\caption{Reducing the scale~$t$, the corresponding diffusion distances become locally unstable, as we move far from the seed point (on the feet). See also Fig.~\ref{fig_OMOTONDO-SIMPLIFIED}.\label{fig:HK-UNDULATIONS}}
\end{figure} 
\subsection{Discretisation of spectral kernels and distances\label{sec:DISCRETE-OPER-DISTANCE}}
Let us consider a (triangular, polygonal, volumetric) mesh \mbox{$\mathcal{M}:=(\mathcal{P},\mathcal{T})$}, which discretises a domain~$\mathcal{N}$, where \mbox{$\mathcal{P}:=\{\mathbf{p}_{i}\}_{i=1}^{n}$} is the set of~$n$ vertices and~$\mathcal{T}$ is the connectivity graph. On~$\mathcal{M}$, a scalar function \mbox{$f:\mathcal{M}\rightarrow\mathbb{R}$} is identified with the vector \mbox{$\mathbf{f}:=(f(\mathbf{p}_{i}))_{i=1}^{n}$} of~$f$-values at~$\mathcal{P}$.

Let us introduce the \emph{Laplacian matrix} \mbox{$\tilde{\mathbf{L}}:=\mathbf{B}^{-1}\mathbf{L}$}, where~$\mathbf{L}$ and~$\mathbf{B}$ are the stiffness and mass matrix respectively (e.g., cotangent~\citep{PINKALL1993}, Voronoi-cotg~\citep{DESBRUN1999}, linear FEM~\citep{reuter:cad06} weights), and the corresponding spectral decomposition is \mbox{$\mathbf{L}\mathbf{X}=\mathbf{B}\mathbf{X}\Lambda$}, \mbox{$\mathbf{X}^{\top}\mathbf{B}\mathbf{X}=\mathbf{I}$}, where \mbox{$\mathbf{X}:=[\mathbf{x}_{1},\ldots,\mathbf{x}_{n}]$} is the eigenvectors' matrix and~$\Lambda$ is the diagonal matrix of the eigenvalues \mbox{$(\lambda_{i})_{i=1}^{n}$}. Analogous discretisations apply to polygonal~\citep{ALEXA2011,HERHOLZ2015} and tetrahedral~\citep{LIAO2009,TONG2003} meshes, or point sets~\citep{LIU2012}.

Under these assumptions, we introduce the spectral kernels and distances (Sect.~\ref{sec:DISCRETE-SPECTRAL-KERNEL}), together with their spectrum-free approximation (Sect.~\ref{sec:DISCRETE-SPECTRAL-APPROXIMATION}), which is based on a rational approximation of the input filter.

\subsubsection{Discrete spectral kernels and distances\label{sec:DISCRETE-SPECTRAL-KERNEL}}
The spectral operator~$\Phi_{1/\rho}$ is discretised by the \emph{spectral kernel matrix}~$\mathbf{K}_{1/\rho}$ such that 
\begin{equation*}
\rho^{-1}(\lambda_{i})
=_{\textrm{Eq. }(\ref{eq:FUNCT-OPER})}\langle\Phi_{1/\rho}\phi_{i},\phi_{j}\rangle_{2}
=\mathbf{x}_{i}^{\top}\mathbf{K}_{1/\rho}^{\top}\mathbf{B}\mathbf{x}_{j},\quad \forall i=1,\ldots,n.
\end{equation*}
Indeed, \mbox{$\mathbf{K}_{1/\rho}=\mathbf{X}\rho^{\dag}(\Lambda)\mathbf{X}^{\top}\mathbf{B}$} is the pseudo-inverse of the \emph{spectral kernel} \mbox{$\mathbf{K}_{\rho}=\rho(\tilde{\mathbf{L}})$}, which is a filtered version of the Laplacian matrix. Here, \mbox{$\rho^{\dag}(\Lambda)=\textrm{diag}(1/\rho(\lambda_{i}))_{i=1}^{n}$} and its entry is null if \mbox{$\rho(\lambda_{i})=0$}. Then, the \emph{discrete spectral distances}~\citep{PATANE-STAR2016,PATANE-CGF2017} are
\begin{equation*}\label{eq:GEN-SPEC-DIST-KERNEL}
d^{2}(\mathbf{p}_{i},\mathbf{p}_{j})
=\|\mathbf{K}_{1/\rho}(\mathbf{e}_{i}-\mathbf{e}_{j})\|_{\mathbf{B}}^{2}
=\sum_{l=1}^{n}\frac{\vert\langle\mathbf{x}_{l},\mathbf{e}_{i}-\mathbf{e}_{j}\rangle_{\mathbf{B}}\vert^{2}}{\rho^{2}(\lambda_{l})},
\end{equation*}
where~$\mathbf{e}_{i}$ is the vector of the canonical basis of~$\mathbb{R}^{n}$, \mbox{$\langle\mathbf{f},\mathbf{g}\rangle_{\mathbf{B}}:=\mathbf{f}^{\top}\mathbf{B}\mathbf{g}$} and \mbox{$\|\mathbf{f}\|_{\mathbf{B}}^{2}:=\mathbf{f}^{\top}\mathbf{B}\mathbf{f}$} are the scalar product and the norm induced by the mass matrix, respectively. In this case, we have derived the spectral distances by applying the continuous expression of the spectral operator instead of the spectral kernel, as done in~\citep{,PATANE-STAR2016,PATANE-CGF2017}.

Analogously to the definitions in Sect.~\ref{sec:SPEC-OP-DIST}, the spectral distance is equal to the norm \mbox{$d(\mathbf{f},\mathbf{g}):=\|\mathbf{u}\|_{\mathbf{B}}$} of the solution to the linear system \mbox{$\mathbf{K}_{\rho}\mathbf{u}=\mathbf{f}-\mathbf{g}$} or to the norm of the vector \mbox{$\mathbf{u}=\mathbf{K}_{1/\rho}(\mathbf{f}-\mathbf{g})$}. In a similar way, the spectral distance between two points reduces to
\begin{equation*}
d(\mathbf{p}_{i},\mathbf{p}_{j})
=\|\mathbf{K}_{1/\varphi}(\mathbf{e}_{i}-\mathbf{e}_{j})\|_{\mathbf{B}}.
\end{equation*}
Indeed, this approximation of the spectral distances involves the spectral kernel only and allows us to bypass numerical inaccuracies due to repeated or close Laplacian eigenvalues~\citep{GOLUB1989} ($\S$~$7$). 
\begin{figure}[t]
\centering
\begin{tabular}{cccc}
\includegraphics[height=75pt]{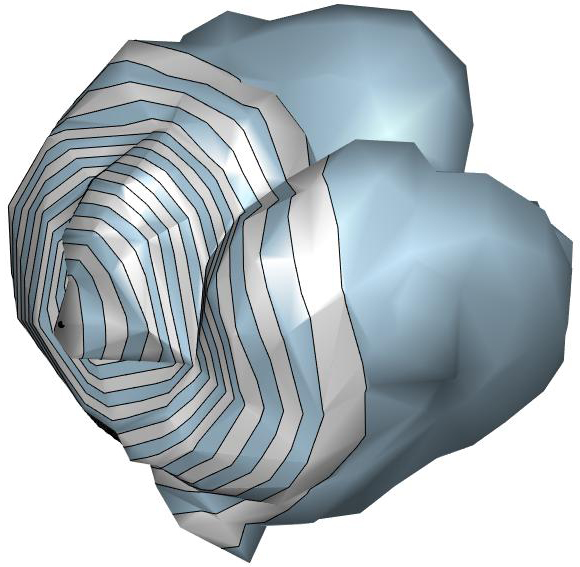}
&\includegraphics[height=75pt]{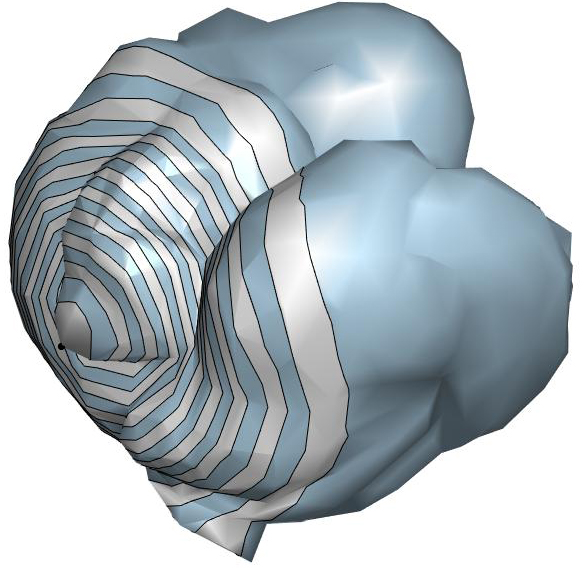}
&\includegraphics[height=75pt]{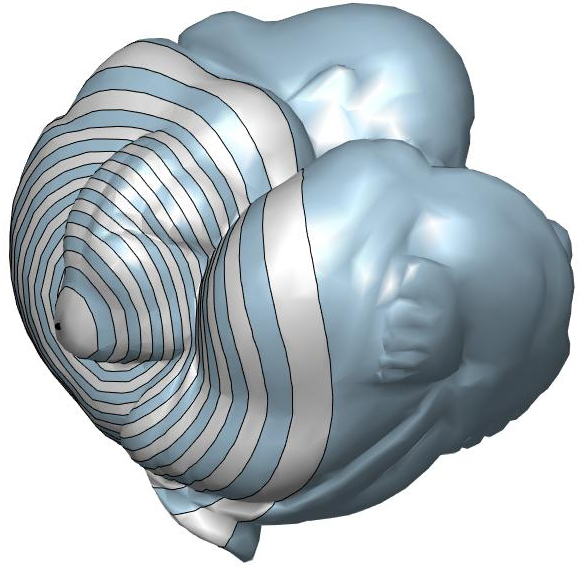}
&\includegraphics[height=75pt]{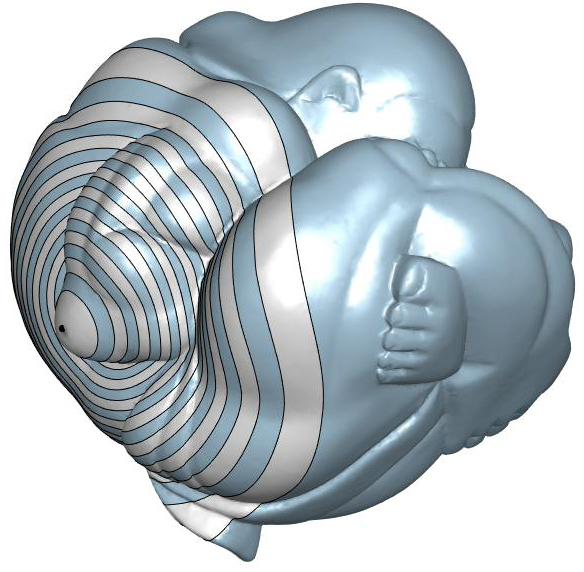}\\
$n=200$ &$n=500$ &$n=2K$ &$n=60K$
\end{tabular}
\caption{Analogous behaviour of the heat kernel with respect to a different resolution ($n$ vertices) of the input shape, in terms of shape and distribution of the level-sets.\label{fig_OMOTONDO-SIMPLIFIED}}
\end{figure}
\subsubsection{Spectral kernel approximation and conditioning\label{sec:DISCRETE-SPECTRAL-APPROXIMATION}}
Since the eigenvalues of the spectral kernel matrix~$\mathbf{K}_{\rho}$ are \mbox{$(\rho(\lambda_{i}))_{i=1}^{n}$}, the approximation of~$\mathbf{K}_{\rho}$ with a new kernel~$\mathbf{K}_{\varphi}$ reduces to the approximation of the corresponding filters with respect to the~$\ell_{\infty}$ norm; in fact,
\begin{equation*}
\|\mathbf{K}_{\rho}-\mathbf{K}_{\varphi}\|_{2}
=\|\mathbf{K}_{\rho-\varphi}\|_{2}
=\max_{i=1,\ldots,n}\{\vert\rho(\lambda_{i})-\varphi(\lambda_{i})\vert\}
\leq\|\rho-\varphi\|_{\infty}.
\end{equation*}
The approximation~$\varphi$ of~$\rho$ is computed on the interval \mbox{$[0,\lambda_{\max}(\tilde{\mathbf{L}})]$}, where the maximum Laplacian eigenvalue is evaluated by the Arnoldi method~\citep{GOLUB1989}, or is set equal to the upper bound \mbox{$\lambda_{\max}(\tilde{\mathbf{L}})\leq \min\{\max_{i}\{\sum_{j}\tilde{L}(i,j)\},\max_{j}\{\sum_{i}\tilde{L}(i,j)\}\}$}~\citep{LEHOUCQ1996,SORENSEN1992}. 

\paragraph*{Conditioning of the spectral kernel}
We now analyse the conditioning of the spectral kernel. Assuming that~$\rho$ is an increasing function (i.e.,~$1/\rho$ is a low pass filter), the conditioning number of the filtered Laplacian matrix is bounded as
\begin{equation*}
\kappa_{2}(\mathbf{K}_{\rho})
=\kappa_{2}(\rho(\tilde{\mathbf{L}}))=\frac{\max_{i=1,\ldots,n}\{\rho(\lambda_{i})\}}{\min_{i=1,\ldots,n}\{\rho(\lambda_{i})\}}=\frac{\|\rho\|_{\infty}}{\rho(0)},
\end{equation*}
and it is ill-conditioned when \mbox{$\rho(0)$} is close to zero or~$\rho$ is unbounded. If~$\rho$ is bounded and \mbox{$\rho(0)$} is not too close to~$0$, then the filtered Laplacian matrix is well-conditioned. If \mbox{$\rho(0)$} is null, then we consider the smallest and not null filtered Laplacian eigenvalue at the denominator of the previous relation. 
\begin{figure*}[t]
\centering
\begin{tabular}{c|c}
$t=10^{-1}$ &~$t=10^{-2}$\\
\hline
\includegraphics[width=165pt]{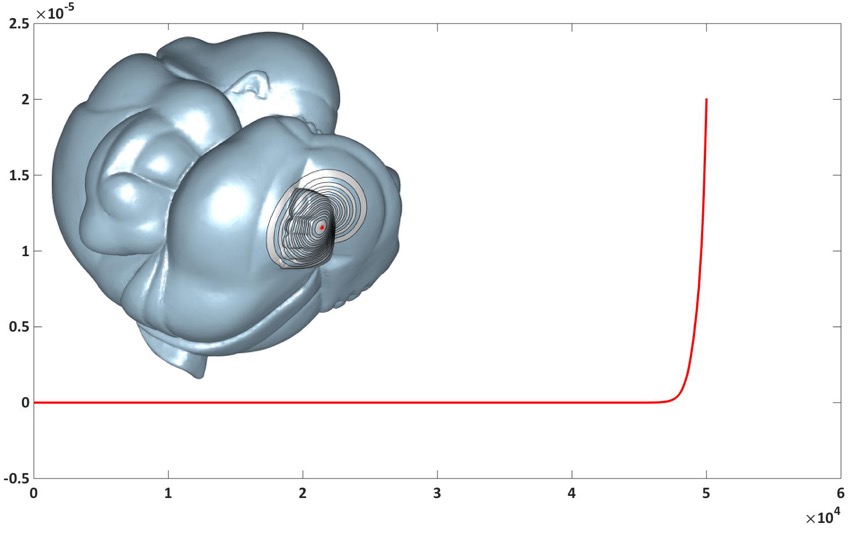}
&\includegraphics[width=165pt]{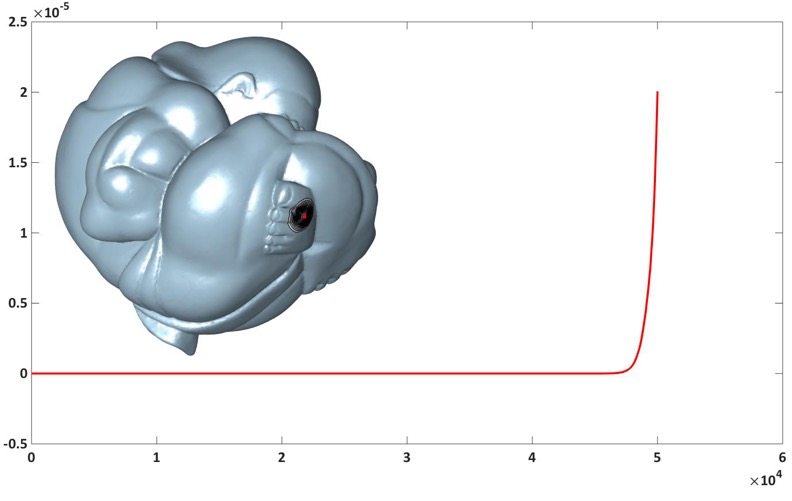}\\
(a) &(b)\\
\includegraphics[width=165pt]{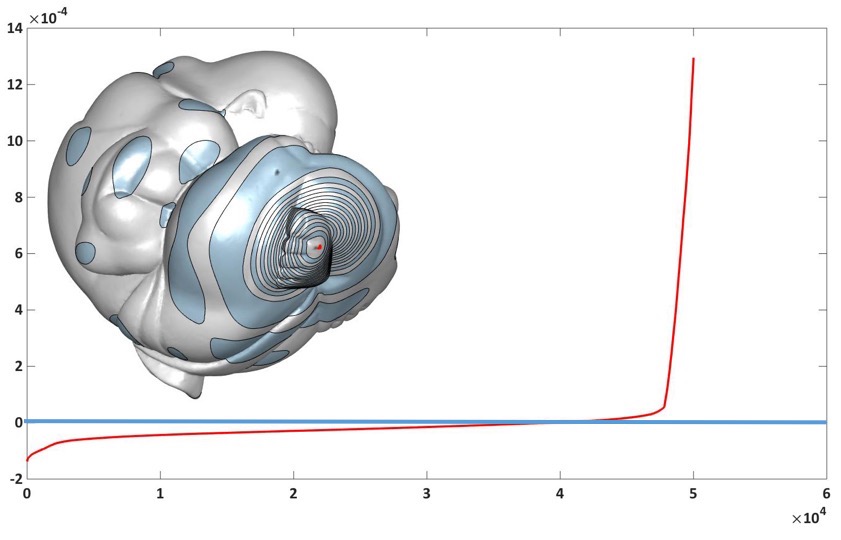}
&\includegraphics[width=165pt]{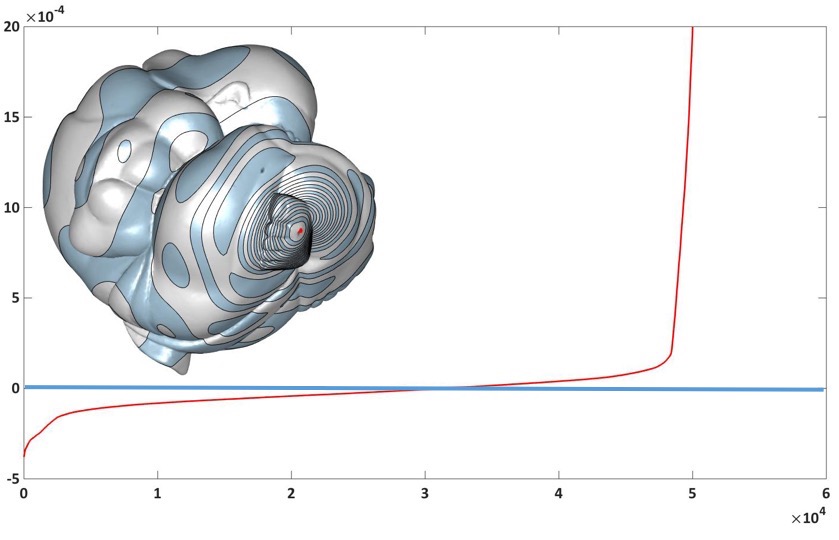}\\
(c) &(d)
\end{tabular}
\caption{(c,d) Small undulations and negative values of the truncated spectral approximation at small scales. These local undulations are not present in the Pad\'e-Chebyshev approximation (a,b), whose values are always positive.\label{fig:GIBBS-PHENOMENON}}
\end{figure*}
\subsection{Spectrum-free approximation of distances: canonical basis\label{sec:UNIFIED-COMPUTATION}}
Analogously to the discussion in Sect.~\ref{sec:RATIONAL-APPROX}, the basis~$\mathcal{B}$ used to represent the rational approximation of the input filter simplifies the computation of the spectral kernel; in fact, we need to compute only the following vectors \mbox{$\tilde{\mathbf{L}}^{i}\mathbf{f}$} and \mbox{$\tilde{(\mathbf{L}}^{i})^{\dag}\mathbf{f}$}. Applying the spectral representation of the Laplacian matrix
\begin{equation*}
\tilde{\mathbf{L}}=\mathbf{X}\Lambda\mathbf{X}^{\top}\mathbf{B},\quad
\Lambda=\textrm{diag}(\lambda_{i})_{i=1}^{n},\quad
\lambda_{1}=0,\quad\lambda_{i}\leq\lambda_{i+1},\quad i=2,\ldots,n-1,
\end{equation*}
we can represent its powers and the corresponding pseudo-inverse matrices as 
\begin{equation*}
\begin{split}
\tilde{\mathbf{L}}^{i}=\mathbf{X}\Lambda^{i}\mathbf{X}^{\top}\mathbf{B}.\quad
(\tilde{\mathbf{L}}^{i})^{\dag}
=\mathbf{X}(\Lambda^{i})^{\dag}\mathbf{X}^{\top}\mathbf{B},\quad
\Lambda^{\dag}=\textrm{diag}(0,\lambda_{i}^{-1})_{i=2}^{n}.
\end{split}
\end{equation*}
To compute \mbox{$\mathbf{g}_{1}=\tilde{\mathbf{L}}^{\dag}\mathbf{f}$}, let us multiply both sides by~$\tilde{\mathbf{L}}$ and notice that
\begin{equation*}
\begin{split}
\tilde{\mathbf{L}}\mathbf{g}_{1}
&=\tilde{\mathbf{L}}\tilde{\mathbf{L}}^{\dag}\mathbf{f}\\
&=\mathbf{X}\Lambda\Lambda^{\dag}\mathbf{X}^{\top}\mathbf{B}\mathbf{f},
\quad\mathbf{X}^{\top}\mathbf{B}\mathbf{X}=\mathbf{I},\\
&=\mathbf{X}(\mathbf{I}-\mathbf{e}_{1}\mathbf{e}_{1}^{\top})\mathbf{X}^{\top}\mathbf{B}\mathbf{f},\quad \mathbf{e}_{1}:=[1,0,\ldots,0]^{\top},\\
&=\mathbf{f}-(\mathbf{1}^{\top}\mathbf{B}\mathbf{f})\mathbf{1},
\quad \mathbf{1}:=[1,1,\ldots,1]^{\top};
\end{split}
\end{equation*}
i.e.,~$\mathbf{g}_{1}$ is the least-squares solution to the sparse and symmetric linear system
\begin{equation*}
\mathbf{L}\mathbf{g}_{1}=\mathbf{B}\mathbf{g}_{0},\quad
\mathbf{g}_{0}:=\mathbf{f}-(\mathbf{1}^{\top}\mathbf{B}\mathbf{f})\mathbf{1},
\end{equation*}
where~$\mathbf{g}_{0}$ is achieved by subtracting to~$\mathbf{f}$ its mean value \mbox{$\langle\mathbf{f},\mathbf{1}\rangle_{\mathbf{B}}$} (c.f., Eq. (\ref{eq:CONTINUOUS-APPROX})). For the general case, we apply the recursive relation
\begin{equation*}
\mathbf{g}_{i}
:=(\tilde{\mathbf{L}}^{\dag})^{i}\mathbf{f}
=\tilde{\mathbf{L}}^{\dag}(\tilde{\mathbf{L}}^{\dag})^{i-1}\mathbf{f}
=\tilde{\mathbf{L}}^{\dag}\mathbf{g}_{i-1},\quad i\geq 2,
\end{equation*}
which reduces to the previous case (Fig.~\ref{fig:BASIS-FUNCTIONS}). In a similar way, \mbox{$\mathbf{g}_{i}=\tilde{\mathbf{L}}^{i}\mathbf{f}$} is calculated by recursively solving the sparse, symmetric, and positive-definite linear systems \mbox{$\mathbf{B}\mathbf{g}_{i+1}=\mathbf{L}\mathbf{g}_{i}$}, \mbox{$i=1,\ldots,r-1$}, with \mbox{$\mathbf{g}_{0}:=\mathbf{f}$}.
\begin{figure}[t]
\centering
\begin{tabular}{cc|cc}
\includegraphics[height=75pt]{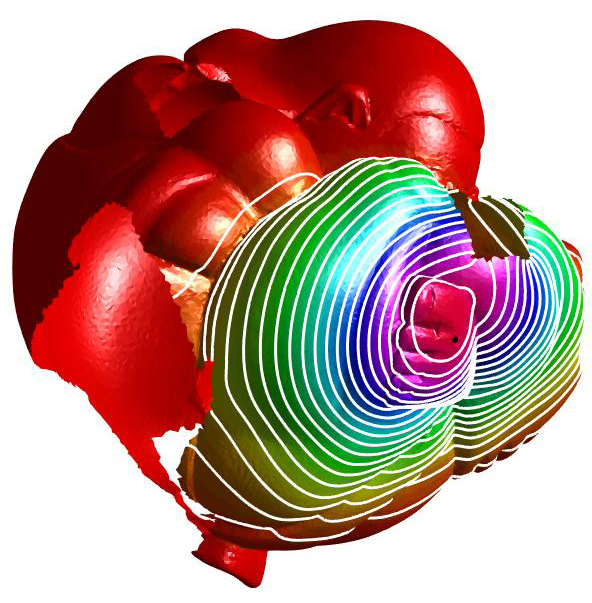}
&\includegraphics[height=75pt]{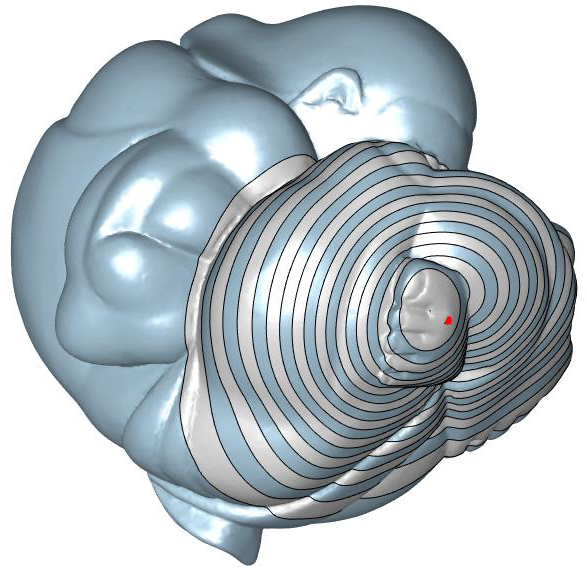}
&\includegraphics[height=75pt]{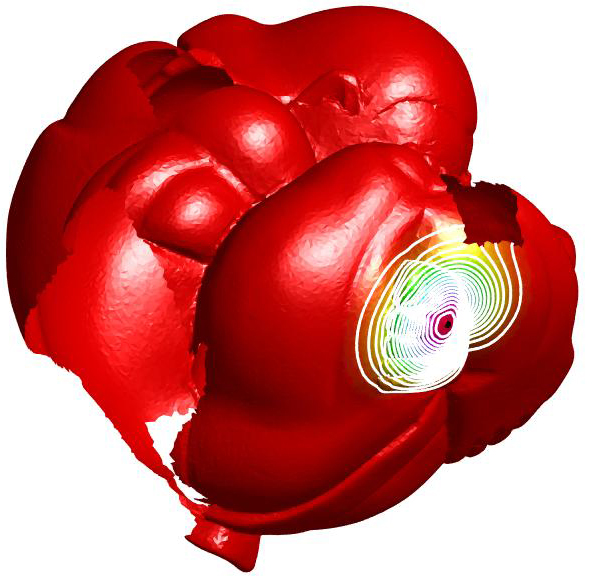}
&\includegraphics[height=75pt]{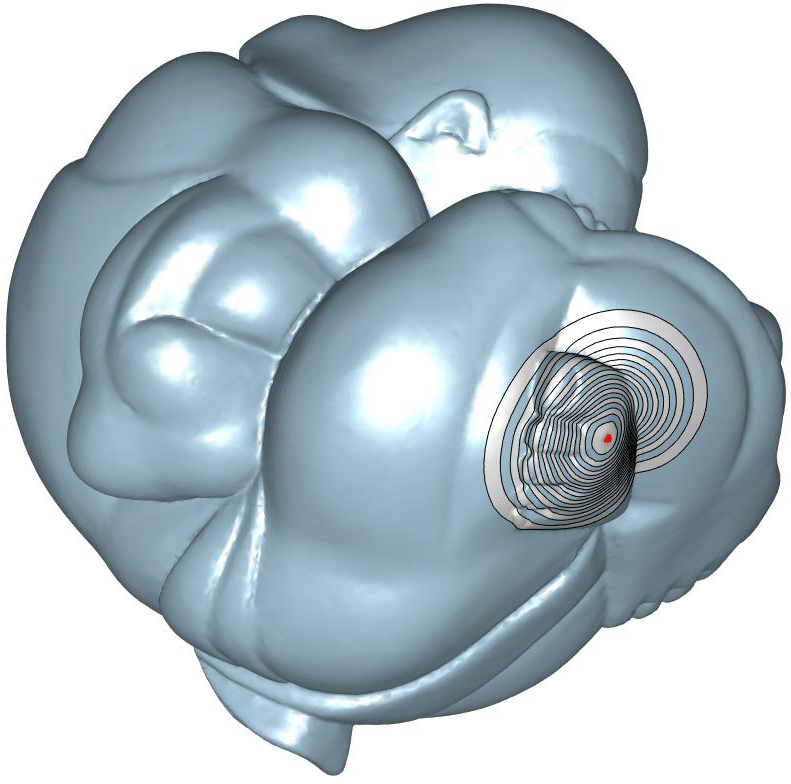}\\
\multicolumn{2}{c}{$t=0.1$} &\multicolumn{2}{c}{$t=0.01$}
\end{tabular}
\caption{Stability of the spectrum-free computation of the diffusion kernel with respect to partially-sampled surfaces.\label{fig:HOLE-ROBUSTNESS}}
\end{figure}
\subsection{Computational cost\label{sec:COMPUTATIONAL-COST}}
The truncated spectral approximation, whose computational cost depends on the sparsity degree of the Laplacian matrix, takes from \mbox{$\mathcal{O}(kn\log n)$} to \mbox{$\mathcal{O}(kn^{2})$} time, where~$k$ is the number of selected eigenpairs. Selecting a rational approximation of the input filter of degree \mbox{$(r,l)$}, the evaluation of the corresponding spectral kernel is reduced to solve~$l$ linear systems whose coefficient matrix is~$\mathbf{B}$ and~$r$ linear systems whose coefficient matrix is~$\mathbf{L}$. Through iterative solvers, the computational cost is \mbox{$\mathcal{O}((r+l)\tau(n))$}, where \mbox{$\tau(n)$} is the cost for the solution of a sparse linear system, which varies from \mbox{$\mathcal{O}(n)$} to \mbox{$\mathcal{O}(n^{2})$}, according to the sparsity of the coefficient matrix, and it is \mbox{$\mathcal{O}(n\log n)$} in the average case.
\begin{figure}[t]
\centering
\begin{tabular}{cc}
(a)\includegraphics[width=120pt]{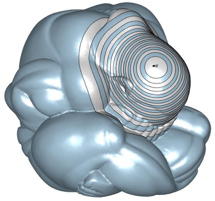}
&(b)\includegraphics[width=150pt]{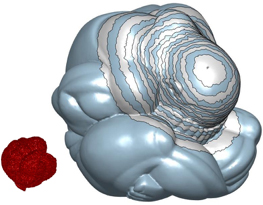}\\
(c)\includegraphics[width=150pt]{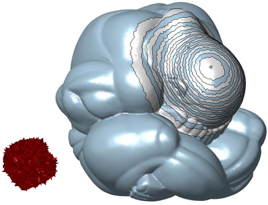}
&(d)\includegraphics[width=150pt]{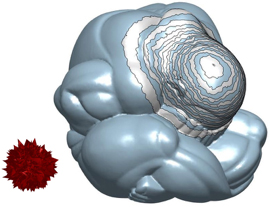}
\end{tabular}
\caption{Behaviour and robustness of the spectrum-free approximation of the diffusion kernel on (b-d) the noisy shapes (bottom part, red surfaces), which have been plotted on (a) the initial shape.\label{fig:SHAPE-NOISE}}
\end{figure}

In a similar way, the computation of the spectral distance between two points has the same order of complexity of the evaluation of the corresponding spectral kernel at the same input points. For the computation of the one-to-all spectral distance, we lump or pre-factorise~$\mathbf{B}$ (if not already diagonal) and pre-factorise~$\mathbf{L}$. Then, the overall computational cost varies from \mbox{$\mathcal{O}(rn)$} ($\mathbf{B}$ diagonal) to \mbox{$\mathcal{O}(n\log n+rn)$} ($\mathbf{B}$ not diagonal), where \mbox{$\mathcal{O}(n\log n)$} is the time for the factorisation of~$\mathbf{L}$. 

In spite of this different computational cost, the truncated spectral approximation (Fig.~\ref{fig:DIFFUSIVE-COMPUTATION}) is affected by small geometric undulations (especially at small scales), the use of heuristics for the selection of the number of Laplacian eigenpairs with respect to the target approximation accuracy, and the scale of features of the input shape. The spectrum-free computation (Fig.~\ref{fig:OMOTONDO-RATIONAL}) generally provides better results in terms of smoothness, regularity, and accuracy of the computed spectral basis. For further comparison examples, we refer the reader to Sect.~\ref{sec:COMPARISON}. 
\begin{figure}[t]
\centering
\begin{tabular}{ccc}
\hline
$t=10^{-2}$ &$t=10^{-1}$ &$t=1$\\
\hline
\includegraphics[height=90pt]{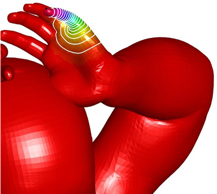}
&\includegraphics[height=90pt]{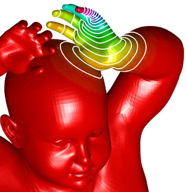}
&\includegraphics[height=90pt]{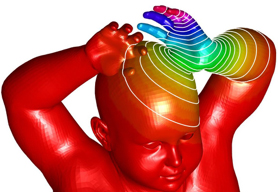}\\
\includegraphics[height=90pt]{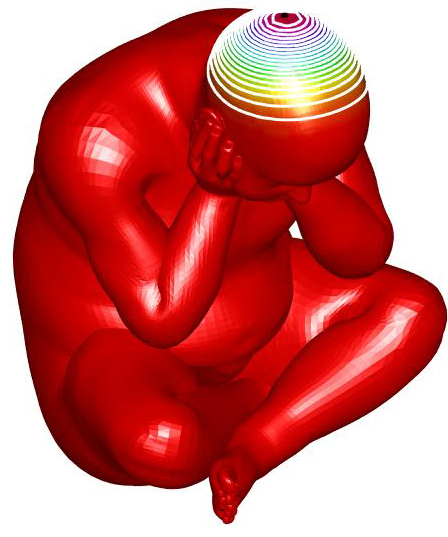}
&\includegraphics[height=90pt]{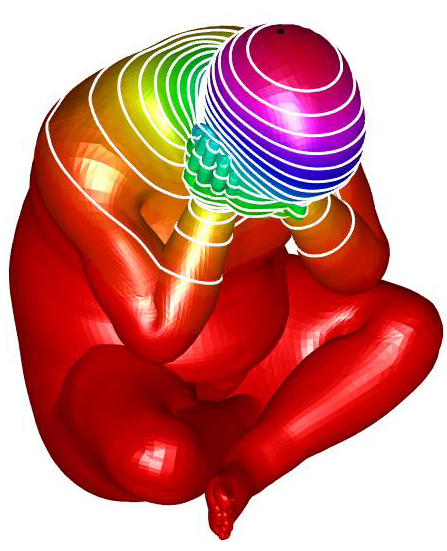}
&\includegraphics[height=90pt]{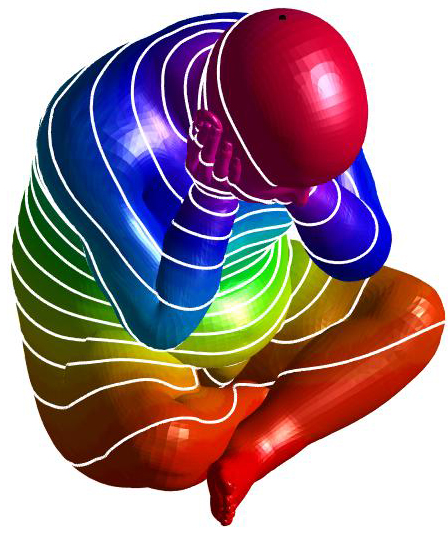}\\
\hline\hline
\end{tabular}
\begin{tabular}{cc}
\includegraphics[height=100pt]{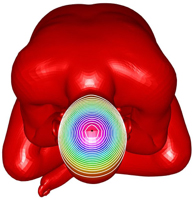}
&\includegraphics[height=100pt]{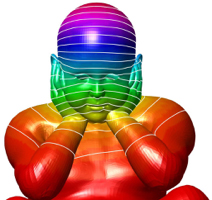}\\
$t=10^{-2}$ &$t=10^{-1}$\\
\end{tabular}
\caption{Robustness of the diffusion kernel at different scales and on self-intersecting surfaces.\label{fig:TOPOLOGICAL-NOISE}}
\end{figure}
\section{Discussion\label{sec:COMPARISON}}
We now discuss the selection of the filter map and main examples on the computation of the Laplacian spectral distances.

The spectral kernels and distances depend only on the behaviour of the filter in the spectral domain and on the Laplacian spectrum: increasing or decreasing the decay of the filter to zero encodes global or local shape details, respectively. Recalling that an arbitrary filter~$1/\rho$ can be approximated in~$\mathcal{R}^{l}_{r}$ with an accuracy of order \mbox{$\mathcal{O}(s^{l+r+1})$} with respect to the~$\ell_{\infty}$-norm (Sect.~\ref{sec:UNIFIED-COMPUTATION}), we can use the space~$\mathcal{R}^{l}_{r}$ to define any filter and the corresponding spectral distances. In this way, we reduce the degree of freedom in the definition of the filter to the selection of \mbox{$(l+r+1)$} coefficients and without losing the richness of the resulting spectral distances. For the selection of the coefficients of the rational filter and the filter frequencies, we can apply the rules proposed in~\citep{KIM2005} for Laplacian spectral smoothing. In particular, the filter frequencies are derived from the dimension of a bounding box placed around the chosen feature~$\mathcal{F}$ and whose axis are aligned with the eigenvectors of the covariance matrix of~$\mathcal{F}$.
\begin{figure}[t]
\centering
\includegraphics[width=320pt]{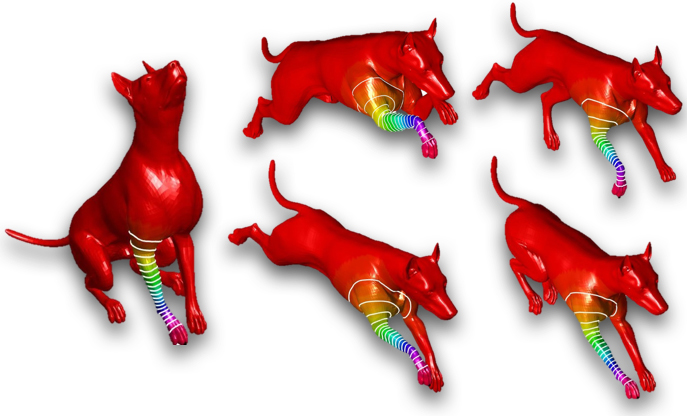}
\caption{Behaviour of the diffusion kernel on (almost) isometric 3D shapes.\label{fig:SHAPE-DEFORMATION}}
\end{figure}

The truncated spectral approximation (c.f., Eq. (\ref{eq:TRUNCATED-SPECTRAL-APPROXIMATION})) of the diffusion, and more generally of the spectral distances, is affected by (i) small undulations (especially at small scales, Fig.~\ref{fig:HK-UNDULATIONS}), (ii) heuristics for the selection of the number of Laplacian eigenpairs with respect to the target approximation accuracy and the scale of the shape features, and (iii) the overall computational cost and storage overhead for the computation of a subpart of the Laplacian spectrum. The comparison of these results with the ones induced by the spectrum-free computation (Fig.~\ref{fig_OMOTONDO-SIMPLIFIED}) shows the improvement of the proposed approach in terms of smoothness, regularity, and accuracy of the computed distances.

At small scales (Fig.~\ref{fig:GIBBS-PHENOMENON}(e,f)), the truncated spectral approximation of the diffusion kernel generally have small negative values as a matter of the slow decay of the exponential filter to zero for small eigenvalues. For the Pad\'e-Chebyshev approximation (Fig.~\ref{fig:GIBBS-PHENOMENON}(a-d)), the kernel values are positive at all the scales, as we approximate the filter with a higher accuracy with respect to the truncated spectral approximation, which applies a low-pass filter. 

Finally, the analogous distribution and shape of the level-sets of the heat kernel and diffusion distance confirm the robustness of the Pad\'e-Chebyshev of the approximation with respect to surface resolution (Fig.~\ref{fig_OMOTONDO-SIMPLIFIED}), partial sampling (Fig.~\ref{fig:HOLE-ROBUSTNESS}), geometric (Fig.~\ref{fig:SHAPE-NOISE}) and topological (Fig.~\ref{fig:TOPOLOGICAL-NOISE}) noise, almost isometric deformations (Fig.~\ref{fig:SHAPE-DEFORMATION}).

\section{Conclusions and future work\label{sec:CONCLUSIONS}}
We have presented a unified approach to the definition, discretisation, and computation of the spectral kernels and distances, which are defined by filtering the Laplacian spectrum and generalise the commute-time, bi-harmonic, diffusion, and wave kernel and distances. Even though analytic filters are easily defined and encode global/local shape properties, the corresponding spectral kernels and distances are not capable of characterising a specific shape class (e.g., humans with respect to four-legs animals, chairs with respect to tables). Indeed, the main topic for future research is a deeper analysis of the constraints on the filter in order to define ``optimal'' spectral kernels and distances for shape comparison. To achieve this goal, we plan to apply and specialise learning methods, as initially investigated in~\citep{AFLALO2011,BOSCAINI2015A}. Another interesting aspect is the study of spectral kernels for the adaptive hash retrieval~\citep{BAI2018}  and the identification of those PDEs, whose Green kernels are useful for shape analysis, as already demonstrated by the bi-harmonic and diffusion kernels.

\paragraph*{Acknowledgements} 
We thank the Reviewers for their thorough review and constructive comments, which helped us to improve the technical part and presentation of the revised paper. This work is supported by the H2020 ERC Advanced Grant CHANGE, grant agreement No. 694515. Shapes are courtesy of the AIM{@}SHAPE Repository, the SHREC 2010 and 2016 data sets; tet-meshes have been generated by the TETGEN software \url{http://wias-berlin.de/software/tetgen/}.

\begin{thebibliography}{54}
\expandafter\ifx\csname natexlab\endcsname\relax\def\natexlab#1{#1}\fi
\providecommand{\url}[1]{\texttt{#1}}
\providecommand{\href}[2]{#2}
\providecommand{\path}[1]{#1}
\providecommand{\DOIprefix}{doi:}
\providecommand{\ArXivprefix}{arXiv:}
\providecommand{\URLprefix}{URL: }
\providecommand{\Pubmedprefix}{pmid:}
\providecommand{\doi}[1]{\href{http://dx.doi.org/#1}{\path{#1}}}
\providecommand{\Pubmed}[1]{\href{pmid:#1}{\path{#1}}}
\providecommand{\bibinfo}[2]{#2}
\ifx\xfnm\relax \def\xfnm[#1]{\unskip,\space#1}\fi
\bibitem[{Aflalo et~al.(2015)Aflalo, Brezis and R.Kimmel}]{AFLALO2015}
\bibinfo{author}{Aflalo, Y.}, \bibinfo{author}{Brezis, H.},
  \bibinfo{author}{R.Kimmel}, \bibinfo{year}{2015}.
\newblock \bibinfo{title}{On the optimality of shape and data representation in
  the spectral domain}.
\newblock \bibinfo{journal}{{SIAM} Journal Imaging Sciences}
  \bibinfo{volume}{8}, \bibinfo{pages}{1141--1160}.
\bibitem[{Aflalo et~al.(2011)Aflalo, Bronstein, Bronstein and
  Kimmel}]{AFLALO2011}
\bibinfo{author}{Aflalo, Y.}, \bibinfo{author}{Bronstein, A.M.},
  \bibinfo{author}{Bronstein, M.M.}, \bibinfo{author}{Kimmel, R.},
  \bibinfo{year}{2011}.
\newblock \bibinfo{title}{Deformable shape retrieval by learning diffusion
  kernels}, in: \bibinfo{booktitle}{Scale space and variational methods in
  Computer Vision}, pp. \bibinfo{pages}{689--700}.
\bibitem[{Alexa and Wardetzky(2011)}]{ALEXA2011}
\bibinfo{author}{Alexa, M.}, \bibinfo{author}{Wardetzky, M.},
  \bibinfo{year}{2011}.
\newblock \bibinfo{title}{Discrete {L}aplacians on general polygonal meshes}.
\newblock \bibinfo{journal}{ACM Trans. on Graphics} \bibinfo{volume}{30}.
\bibitem[{Bahonar et~al.(2018)Bahonar, Mirzaei and Wilson}]{BAHONAR2018}
\bibinfo{author}{Bahonar, H.}, \bibinfo{author}{Mirzaei, A.},
  \bibinfo{author}{Wilson, R.C.}, \bibinfo{year}{2018}.
\newblock \bibinfo{title}{Diffusion wavelet embedding: a multi-resolution
  approach for graph embedding in vector space}.
\newblock \bibinfo{journal}{Pattern Recognition} \bibinfo{volume}{74},
  \bibinfo{pages}{518 -- 530}.
\bibitem[{Bai et~al.(2018)Bai, Yan, Yang, Bai, Zhou and Hancock}]{BAI2018}
\bibinfo{author}{Bai, X.}, \bibinfo{author}{Yan, C.}, \bibinfo{author}{Yang,
  H.}, \bibinfo{author}{Bai, L.}, \bibinfo{author}{Zhou, J.},
  \bibinfo{author}{Hancock, E.R.}, \bibinfo{year}{2018}.
\newblock \bibinfo{title}{Adaptive hash retrieval with kernel-based
  similarity}.
\newblock \bibinfo{journal}{Pattern Recognition} \bibinfo{volume}{75},
  \bibinfo{pages}{136 -- 148}.
\newblock \bibinfo{note}{Distance Metric Learning for Pattern Recognition}.
\bibitem[{Belkin and Niyogi(2003)}]{BELKIN2003}
\bibinfo{author}{Belkin, M.}, \bibinfo{author}{Niyogi, P.},
  \bibinfo{year}{2003}.
\newblock \bibinfo{title}{{L}aplacian eigenmaps for dimensionality reduction
  and data representation}.
\newblock \bibinfo{journal}{Neural Computations} \bibinfo{volume}{15},
  \bibinfo{pages}{1373--1396}.
\bibitem[{Berard et~al.(1994)Berard, Besson and Gallot}]{BERARD1984}
\bibinfo{author}{Berard, P.}, \bibinfo{author}{Besson, G.},
  \bibinfo{author}{Gallot, S.}, \bibinfo{year}{1994}.
\newblock \bibinfo{title}{Embedding {R}iemannian manifolds by their heat
  kernel}.
\newblock \bibinfo{journal}{Geometric and Functional Analysis}
  \bibinfo{volume}{4}, \bibinfo{pages}{373--398}.
\bibitem[{Boscaini et~al.(2015)Boscaini, Masci, Melzi, Bronstein, Castellani
  and Vandergheynst}]{BOSCAINI2015A}
\bibinfo{author}{Boscaini, D.}, \bibinfo{author}{Masci, J.},
  \bibinfo{author}{Melzi, S.}, \bibinfo{author}{Bronstein, M.M.},
  \bibinfo{author}{Castellani, U.}, \bibinfo{author}{Vandergheynst, P.},
  \bibinfo{year}{2015}.
\newblock \bibinfo{title}{Learning class-specific descriptors for deformable
  shapes using localized spectral convolutional networks}.
\newblock \bibinfo{journal}{Computer Graphics Forum} \bibinfo{volume}{34},
  \bibinfo{pages}{13--23}.
\bibitem[{Boscaini et~al.(2016)Boscaini, Masci, Rodol\`a, Bronstein and
  Cremers}]{BOSCAINI2016}
\bibinfo{author}{Boscaini, D.}, \bibinfo{author}{Masci, J.},
  \bibinfo{author}{Rodol\`a, E.}, \bibinfo{author}{Bronstein, M.M.},
  \bibinfo{author}{Cremers, D.}, \bibinfo{year}{2016}.
\newblock \bibinfo{title}{Anisotropic diffusion descriptors}.
\newblock \bibinfo{journal}{Computer Graphics Forum} \bibinfo{volume}{35},
  \bibinfo{pages}{431--441}.
\bibitem[{Bronstein et~al.(2010)Bronstein, Bronstein, Kimmel, Mahmoudi and
  Sapiro}]{BRONSTEIN2009}
\bibinfo{author}{Bronstein, A.}, \bibinfo{author}{Bronstein, M.},
  \bibinfo{author}{Kimmel, R.}, \bibinfo{author}{Mahmoudi, M.},
  \bibinfo{author}{Sapiro, G.}, \bibinfo{year}{2010}.
\newblock \bibinfo{title}{A {G}romov-{H}ausdorff framework with diffusion
  geometry for topologically-robust non-rigid shape matching}.
\newblock \bibinfo{journal}{Intern. Journal of Computer Vision}
  \bibinfo{volume}{2-3}, \bibinfo{pages}{266--286}.
\bibitem[{Bronstein et~al.(2011)Bronstein, Bronstein, Ovsjanikov and
  Guibas}]{BRONSTEIN2010-TOG}
\bibinfo{author}{Bronstein, A.M.}, \bibinfo{author}{Bronstein, M.M.},
  \bibinfo{author}{Ovsjanikov, M.}, \bibinfo{author}{Guibas, L.J.},
  \bibinfo{year}{2011}.
\newblock \bibinfo{title}{Shape {G}oogle: geometric words and expressions for
  invariant shape retrieval}.
\newblock \bibinfo{journal}{ACM Trans. on Graphics} \bibinfo{volume}{30}.
\bibitem[{Bronstein and Bronstein(2011)}]{BRONSTEIN-PAMI2011}
\bibinfo{author}{Bronstein, M.}, \bibinfo{author}{Bronstein, A.},
  \bibinfo{year}{2011}.
\newblock \bibinfo{title}{Shape recognition with spectral distances}.
\newblock \bibinfo{journal}{IEEE Trans. on Pattern Analysis and Machine
  Intelligence} \bibinfo{volume}{33}, \bibinfo{pages}{1065 --1071}.
\bibitem[{Cosmo et~al.(2017)Cosmo, Rodol\'{a}, Albarelli, M{\'e}moli and
  Cremers}]{COSMO2016}
\bibinfo{author}{Cosmo, L.}, \bibinfo{author}{Rodol\'{a}, E.},
  \bibinfo{author}{Albarelli, A.}, \bibinfo{author}{M{\'e}moli, F.},
  \bibinfo{author}{Cremers, D.}, \bibinfo{year}{2017}.
\newblock \bibinfo{title}{Consistent partial matching of shape collections via
  sparse modeling}.
\newblock \bibinfo{journal}{Computer Graphics Forum} \bibinfo{volume}{36},
  \bibinfo{pages}{209--221}.
\bibitem[{Crane et~al.(2013)Crane, Weischedel and Wardetzky}]{CRANE2013}
\bibinfo{author}{Crane, K.}, \bibinfo{author}{Weischedel, C.},
  \bibinfo{author}{Wardetzky, M.}, \bibinfo{year}{2013}.
\newblock \bibinfo{title}{Geodesics in heat: A new approach to computing
  distance based on heat flow}.
\newblock \bibinfo{journal}{ACM Trans. on Graphics} \bibinfo{volume}{32},
  \bibinfo{pages}{152:1--152:11}.
\bibitem[{Desbrun et~al.(1999)Desbrun, Meyer, Schr\"{o}der and
  Barr}]{DESBRUN1999}
\bibinfo{author}{Desbrun, M.}, \bibinfo{author}{Meyer, M.},
  \bibinfo{author}{Schr\"{o}der, P.}, \bibinfo{author}{Barr, A.H.},
  \bibinfo{year}{1999}.
\newblock \bibinfo{title}{Implicit fairing of irregular meshes using diffusion
  and curvature flow}, in: \bibinfo{booktitle}{ACM Siggraph}, pp.
  \bibinfo{pages}{317--324}.
\bibitem[{ElGhawalby and Hancock(2015)}]{ELGHAWALBY2015}
\bibinfo{author}{ElGhawalby, H.}, \bibinfo{author}{Hancock, E.R.},
  \bibinfo{year}{2015}.
\newblock \bibinfo{title}{Heat kernel embeddings, differential geometry and
  graph structure}.
\newblock \bibinfo{journal}{Axioms} \bibinfo{volume}{4},
  \bibinfo{pages}{275--293}.
\bibitem[{Gebal et~al.(2009)Gebal, B{\ae}rentzen, Aan{\ae}s and
  Larsen}]{GEBAL2009}
\bibinfo{author}{Gebal, K.}, \bibinfo{author}{B{\ae}rentzen, J.A.},
  \bibinfo{author}{Aan{\ae}s, H.}, \bibinfo{author}{Larsen, R.},
  \bibinfo{year}{2009}.
\newblock \bibinfo{title}{Shape analysis using the auto diffusion function}.
\newblock \bibinfo{journal}{Computer Graphics Forum} \bibinfo{volume}{28},
  \bibinfo{pages}{1405--1413}.
\bibitem[{Ghosh et~al.(2018)Ghosh, Das, Gonzalves, Quaresma and
  Kundu}]{GHOSH2018}
\bibinfo{author}{Ghosh, S.}, \bibinfo{author}{Das, N.},
  \bibinfo{author}{Gonzalves, T.}, \bibinfo{author}{Quaresma, P.},
  \bibinfo{author}{Kundu, M.}, \bibinfo{year}{2018}.
\newblock \bibinfo{title}{The journey of graph kernels through two decades}.
\newblock \bibinfo{journal}{Computer Science Review} \bibinfo{volume}{27},
  \bibinfo{pages}{88 -- 111}.
\bibitem[{de~Goes et~al.(2008)de~Goes, Goldenstein and Velho}]{GOES2008}
\bibinfo{author}{de~Goes, F.}, \bibinfo{author}{Goldenstein, S.},
  \bibinfo{author}{Velho, L.}, \bibinfo{year}{2008}.
\newblock \bibinfo{title}{A hierarchical segmentation of articulated bodies}.
\newblock \bibinfo{journal}{Computer Graphics Forum} \bibinfo{volume}{27},
  \bibinfo{pages}{1349--1356}.
\bibitem[{Golub and VanLoan(1989)}]{GOLUB1989}
\bibinfo{author}{Golub, G.}, \bibinfo{author}{VanLoan, G.},
  \bibinfo{year}{1989}.
\newblock \bibinfo{title}{Matrix Computations}.
\newblock \bibinfo{publisher}{John Hopkins University Press, 2nd Edition}.
\bibitem[{Hammond et~al.(2011)Hammond, Vandergheynst and
  Gribonval}]{HAMMOND2010}
\bibinfo{author}{Hammond, D.K.}, \bibinfo{author}{Vandergheynst, P.},
  \bibinfo{author}{Gribonval, R.}, \bibinfo{year}{2011}.
\newblock \bibinfo{title}{Wavelets on graphs via spectral graph theory}.
\newblock \bibinfo{journal}{Applied and Computational Harmonic Analysis}
  \bibinfo{volume}{30}, \bibinfo{pages}{129 -- 150}.
\bibitem[{Herholz et~al.(2015)Herholz, Kyprianidis and Alexa}]{HERHOLZ2015}
\bibinfo{author}{Herholz, P.}, \bibinfo{author}{Kyprianidis, J.E.},
  \bibinfo{author}{Alexa, M.}, \bibinfo{year}{2015}.
\newblock \bibinfo{title}{Perfect {L}aplacians for polygon meshes}.
\newblock \bibinfo{journal}{Computer Graphics Forum} \bibinfo{volume}{34},
  \bibinfo{pages}{211--218}.
\bibitem[{Hou and Qin(2012)}]{HOU2012}
\bibinfo{author}{Hou, T.}, \bibinfo{author}{Qin, H.}, \bibinfo{year}{2012}.
\newblock \bibinfo{title}{Continuous and discrete {M}exican hat wavelet
  transforms on manifolds}.
\newblock \bibinfo{journal}{Graphical Models} \bibinfo{volume}{74},
  \bibinfo{pages}{221--232}.
\bibitem[{Kim and Rossignac(2005)}]{KIM2005}
\bibinfo{author}{Kim, B.}, \bibinfo{author}{Rossignac, J.},
  \bibinfo{year}{2005}.
\newblock \bibinfo{title}{Geofilter: Geometric selection of mesh filter
  parameters}.
\newblock \bibinfo{journal}{Computer Graphics Forum} \bibinfo{volume}{24},
  \bibinfo{pages}{295--302}.
\bibitem[{Lafon et~al.(2006)Lafon, Keller and Coifman}]{LAFON2006}
\bibinfo{author}{Lafon, S.}, \bibinfo{author}{Keller, Y.},
  \bibinfo{author}{Coifman, R.R.}, \bibinfo{year}{2006}.
\newblock \bibinfo{title}{Data fusion and multicue data matching by diffusion
  maps}.
\newblock \bibinfo{journal}{IEEE Trans. on Pattern Analysis Machine
  Intelligence} \bibinfo{volume}{28}, \bibinfo{pages}{1784--1797}.
\bibitem[{Lehoucq and Sorensen(1996)}]{LEHOUCQ1996}
\bibinfo{author}{Lehoucq, R.}, \bibinfo{author}{Sorensen, D.C.},
  \bibinfo{year}{1996}.
\newblock \bibinfo{title}{Deflation techniques for an implicitly re-started
  {A}rnoldi iteration}.
\newblock \bibinfo{journal}{SIAM Journal of Matrix Analysis and Applications}
  \bibinfo{volume}{17}, \bibinfo{pages}{789--821}.
\bibitem[{Liao et~al.(2009)Liao, Tong, Dong and Zhu}]{LIAO2009}
\bibinfo{author}{Liao, S.}, \bibinfo{author}{Tong, R.}, \bibinfo{author}{Dong,
  J.}, \bibinfo{author}{Zhu, F.}, \bibinfo{year}{2009}.
\newblock \bibinfo{title}{Gradient field based inhomogeneous volumetric mesh
  deformation for maxillofacial surgery simulation}.
\newblock \bibinfo{journal}{Computers $\&$ Graphics} \bibinfo{volume}{33},
  \bibinfo{pages}{424 -- 432}.
\bibitem[{Lipman et~al.(2010)Lipman, Rustamov and Funkhouser}]{LIPMAN2010}
\bibinfo{author}{Lipman, Y.}, \bibinfo{author}{Rustamov, R.M.},
  \bibinfo{author}{Funkhouser, T.A.}, \bibinfo{year}{2010}.
\newblock \bibinfo{title}{Biharmonic distance}.
\newblock \bibinfo{journal}{{ACM} Trans. on Graphics} \bibinfo{volume}{29},
  \bibinfo{pages}{27:1--27:11}.
\bibitem[{Liu et~al.(2012)Liu, Prabhakaran and Guo}]{LIU2012}
\bibinfo{author}{Liu, Y.}, \bibinfo{author}{Prabhakaran, B.},
  \bibinfo{author}{Guo, X.}, \bibinfo{year}{2012}.
\newblock \bibinfo{title}{Point-based manifold harmonics}.
\newblock \bibinfo{journal}{IEEE Trans. on Visualization and Computer Graphics}
  \bibinfo{volume}{18}, \bibinfo{pages}{1693 --1703}.
\bibitem[{Luo et~al.(2003)Luo, Wilson and Hancock}]{LUO2003}
\bibinfo{author}{Luo, B.}, \bibinfo{author}{Wilson, R.C.},
  \bibinfo{author}{Hancock, E.R.}, \bibinfo{year}{2003}.
\newblock \bibinfo{title}{Spectral embedding of graphs}.
\newblock \bibinfo{journal}{Pattern Recognition} \bibinfo{volume}{36},
  \bibinfo{pages}{2213--2230}.
\bibitem[{Mahmoudi and Sapiro(2009)}]{MAHMOUDI2009}
\bibinfo{author}{Mahmoudi, M.}, \bibinfo{author}{Sapiro, G.},
  \bibinfo{year}{2009}.
\newblock \bibinfo{title}{Three-dimensional point cloud recognition via
  distributions of geometric distances}.
\newblock \bibinfo{journal}{Graphical Models} \bibinfo{volume}{71},
  \bibinfo{pages}{22--31}.
\bibitem[{Memoli(2011)}]{MEMOLI2011}
\bibinfo{author}{Memoli, F.}, \bibinfo{year}{2011}.
\newblock \bibinfo{title}{A spectral notion of {G}romov-{W}asserstein distance
  and related methods}.
\newblock \bibinfo{journal}{Applied and Computational Harmonic Analysis}
  \bibinfo{volume}{30}, \bibinfo{pages}{363 -- 401}.
\bibitem[{M\`emoli and Sapiro(2005)}]{Memoli05}
\bibinfo{author}{M\`emoli, F.}, \bibinfo{author}{Sapiro, G.},
  \bibinfo{year}{2005}.
\newblock \bibinfo{title}{A theoretical and computational framework for
  isometry invariant recognition of point cloud data}.
\newblock \bibinfo{journal}{Foundations of Computational Mathematics}
  \bibinfo{volume}{5}, \bibinfo{pages}{313--347}.
\bibitem[{Ng et~al.(2001)Ng, Jordan and Weiss}]{NG2001}
\bibinfo{author}{Ng, A.Y.}, \bibinfo{author}{Jordan, M.I.},
  \bibinfo{author}{Weiss, Y.}, \bibinfo{year}{2001}.
\newblock \bibinfo{title}{On spectral clustering: analysis and an algorithm},
  in: \bibinfo{booktitle}{Advances in Neural Information Processing Systems
  14}, \bibinfo{publisher}{MIT Press}. pp. \bibinfo{pages}{849--856}.
\bibitem[{Ovsjanikov et~al.(2010)Ovsjanikov, M\`erigot, M\`emoli and
  Guibas}]{OVSJANKOV2010}
\bibinfo{author}{Ovsjanikov, M.}, \bibinfo{author}{M\`erigot, Q.},
  \bibinfo{author}{M\`emoli, F.}, \bibinfo{author}{Guibas, L.},
  \bibinfo{year}{2010}.
\newblock \bibinfo{title}{One point isometric matching with the heat kernel}.
\newblock \bibinfo{journal}{Computer Graphics Forum} \bibinfo{volume}{29},
  \bibinfo{pages}{1555--1564}.
\bibitem[{Patan\`e(2013)}]{PATANE-CAGD2013}
\bibinfo{author}{Patan\`e, G.}, \bibinfo{year}{2013}.
\newblock \bibinfo{title}{{wFEM} heat kernel: discretization and applications
  to shape analysis and retrieval}.
\newblock \bibinfo{journal}{Computer-Aided Geometric Design}
  \bibinfo{volume}{30}, \bibinfo{pages}{276--295}.
\bibitem[{Patan{\`{e}}(2014)}]{PATANE2014-PRL}
\bibinfo{author}{Patan{\`{e}}, G.}, \bibinfo{year}{2014}.
\newblock \bibinfo{title}{{L}aplacian spectral distances and kernels on
  {{{{3D}}}} shapes}.
\newblock \bibinfo{journal}{Pattern Recognition Letters} \bibinfo{volume}{47},
  \bibinfo{pages}{102--110}.
\bibitem[{Patan{\`{e}}(2016)}]{PATANE-STAR2016}
\bibinfo{author}{Patan{\`{e}}, G.}, \bibinfo{year}{2016}.
\newblock \bibinfo{title}{{STAR} - {L}aplacian spectral kernels and distances
  for geometry processing and shape analysis}.
\newblock \bibinfo{journal}{Computer Graphics Forum} \bibinfo{volume}{35},
  \bibinfo{pages}{599--624}.
\bibitem[{Patan{\`{e}}(2017)}]{PATANE-CGF2017}
\bibinfo{author}{Patan{\`{e}}, G.}, \bibinfo{year}{2017}.
\newblock \bibinfo{title}{Accurate and efficient computation of {L}aplacian
  spectral distances and kernels}.
\newblock \bibinfo{journal}{Computer Graphics Forum} \bibinfo{volume}{36},
  \bibinfo{pages}{184--196}.
\bibitem[{Pinkall and Polthier(1993)}]{PINKALL1993}
\bibinfo{author}{Pinkall, U.}, \bibinfo{author}{Polthier, K.},
  \bibinfo{year}{1993}.
\newblock \bibinfo{title}{Computing discrete minimal surfaces and their
  conjugates}.
\newblock \bibinfo{journal}{Experimental Mathematics} \bibinfo{volume}{2},
  \bibinfo{pages}{15--36}.
\bibitem[{Ramani and Sinha(2013)}]{SINHA2013}
\bibinfo{author}{Ramani, K.}, \bibinfo{author}{Sinha, A.},
  \bibinfo{year}{2013}.
\newblock \bibinfo{title}{Multiscale kernels using random walks}.
\newblock \bibinfo{journal}{Computer Graphics Forum} \bibinfo{volume}{33},
  \bibinfo{pages}{164--177}.
\bibitem[{Reuter et~al.(2006)Reuter, Wolter and Peinecke}]{reuter:cad06}
\bibinfo{author}{Reuter, M.}, \bibinfo{author}{Wolter, F.E.},
  \bibinfo{author}{Peinecke, N.}, \bibinfo{year}{2006}.
\newblock \bibinfo{title}{{L}aplace-{B}eltrami spectra as {S}hape-{DNA} of
  surfaces and solids}.
\newblock \bibinfo{journal}{Computer-Aided Design} \bibinfo{volume}{38},
  \bibinfo{pages}{342--366}.
\bibitem[{Rodol\`a et~al.(2017)Rodol\`a, Cosmo, Bronstein, Torsello and
  Cremers}]{RODOLA2016}
\bibinfo{author}{Rodol\`a, E.}, \bibinfo{author}{Cosmo, L.},
  \bibinfo{author}{Bronstein, M.M.}, \bibinfo{author}{Torsello, A.},
  \bibinfo{author}{Cremers, D.}, \bibinfo{year}{2017}.
\newblock \bibinfo{title}{Partial functional correspondence}.
\newblock \bibinfo{journal}{Computer Graphics Forum} \bibinfo{volume}{36},
  \bibinfo{pages}{222--236}.
\bibitem[{Rosenberg(1997)}]{ROSENBERG1997}
\bibinfo{author}{Rosenberg, S.}, \bibinfo{year}{1997}.
\newblock \bibinfo{title}{The {L}aplacian on a Riemannian Manifold}.
\newblock \bibinfo{publisher}{Cambridge University Press}.
\bibitem[{Rustamov(2011)}]{RUSTAMOV2011}
\bibinfo{author}{Rustamov, R.M.}, \bibinfo{year}{2011}.
\newblock \bibinfo{title}{Multiscale biharmonic kernels}.
\newblock \bibinfo{journal}{Computer Graphics Forum} \bibinfo{volume}{30},
  \bibinfo{pages}{1521--1531}.
\bibitem[{Shi and Malik(2000)}]{SHI2000}
\bibinfo{author}{Shi, J.}, \bibinfo{author}{Malik, J.}, \bibinfo{year}{2000}.
\newblock \bibinfo{title}{Normalized cuts and image segmentation}.
\newblock \bibinfo{journal}{IEEE Trans. on Pattern Analysis and Machine
  Intelligence} \bibinfo{volume}{22}, \bibinfo{pages}{888--905}.
\bibitem[{Singer(2006)}]{SINGER2006}
\bibinfo{author}{Singer, A.}, \bibinfo{year}{2006}.
\newblock \bibinfo{title}{From graph to manifold {L}aplacian: the convergence
  rate}.
\newblock \bibinfo{journal}{Applied and Computational Harmonic Analysis}
  \bibinfo{volume}{21}, \bibinfo{pages}{128 -- 134}.
\bibitem[{Sorensen(1992)}]{SORENSEN1992}
\bibinfo{author}{Sorensen, D.C.}, \bibinfo{year}{1992}.
\newblock \bibinfo{title}{Implicit application of polynomial filters in a
  k-step arnoldi method}.
\newblock \bibinfo{journal}{SIAM Journal of Matrix Analysis and Applications}
  \bibinfo{volume}{13}, \bibinfo{pages}{357--385}.
\bibitem[{Spielman and Teng(2007)}]{SPIELMAN1996}
\bibinfo{author}{Spielman, D.A.}, \bibinfo{author}{Teng, S.H.},
  \bibinfo{year}{2007}.
\newblock \bibinfo{title}{Spectral partitioning works: planar graphs and finite
  element meshes}.
\newblock \bibinfo{journal}{Linear Algebra and its Applications}
  \bibinfo{volume}{421}, \bibinfo{pages}{284--305}.
\bibitem[{Sun et~al.(2009)Sun, Ovsjanikov and Guibas}]{SUN2009}
\bibinfo{author}{Sun, J.}, \bibinfo{author}{Ovsjanikov, M.},
  \bibinfo{author}{Guibas, L.J.}, \bibinfo{year}{2009}.
\newblock \bibinfo{title}{A concise and provably informative multi-scale
  signature based on heat diffusion}.
\newblock \bibinfo{journal}{Computer Graphics Forum} \bibinfo{volume}{28},
  \bibinfo{pages}{1383--1392}.
\bibitem[{Tong et~al.(2003)Tong, Lombeyda, Hirani and Desbrun}]{TONG2003}
\bibinfo{author}{Tong, Y.}, \bibinfo{author}{Lombeyda, S.},
  \bibinfo{author}{Hirani, A.N.}, \bibinfo{author}{Desbrun, M.},
  \bibinfo{year}{2003}.
\newblock \bibinfo{title}{Discrete multiscale vector field decomposition}.
\newblock \bibinfo{journal}{ACM Trans. on Graphics} \bibinfo{volume}{22},
  \bibinfo{pages}{445--452}.
\bibitem[{Varga(1990)}]{VARGA1990}
\bibinfo{author}{Varga, R.}, \bibinfo{year}{1990}.
\newblock \bibinfo{title}{Scientific computation on mathematical problems and
  conjectures}.
\newblock \bibinfo{publisher}{SIAM, CBMS-NSF regional conference series in
  applied mathematics}.
\bibitem[{Xiao et~al.(2010)Xiao, Hancock and Wilsonb}]{XIAO2010}
\bibinfo{author}{Xiao, B.}, \bibinfo{author}{Hancock, E.R.},
  \bibinfo{author}{Wilsonb, R.}, \bibinfo{year}{2010}.
\newblock \bibinfo{title}{Geometric characterization and clustering of graphs
  using heat kernel embeddings}.
\newblock \bibinfo{journal}{Image and Vision Computing} \bibinfo{volume}{28},
  \bibinfo{pages}{1003 -- 1021}.
\bibitem[{Zhu et~al.(2003)Zhu, Ghahramani and Lafferty}]{ZHU2003}
\bibinfo{author}{Zhu, X.}, \bibinfo{author}{Ghahramani, Z.},
  \bibinfo{author}{Lafferty, J.}, \bibinfo{year}{2003}.
\newblock \bibinfo{title}{Semi-supervised learning using gaussian fields and
  harmonic functions}, in: \bibinfo{booktitle}{Intern. Conf. on Machine
  Learning}, pp. \bibinfo{pages}{912--919}.

\end{thebibliography}
\end{document}